\documentclass[12pt]{article}
\usepackage{amssymb,amsmath,latexsym}
\allowdisplaybreaks

\begin{document}

\begin{doublespace}
\def\1{{\bf 1}}
\def\ind{{\bf 1}}
\def\nn{\nonumber}
\newcommand{\I}{\mathbf{1}}

\def\sA {{\cal A}} \def\sB {{\cal B}} \def\sC {{\cal C}}
\def\sD {{\cal D}} \def\sE {{\cal E}} \def\sF {{\cal F}}
\def\sG {{\cal G}} \def\sH {{\cal H}} \def\sI {{\cal I}}
\def\sJ {{\cal J}} \def\sK {{\cal K}} \def\sL {{\cal L}}
\def\sM {{\cal M}} \def\sN {{\cal N}} \def\sO {{\cal O}}
\def\sP {{\cal P}} \def\sQ {{\cal Q}} \def\sR {{\cal R}}
\def\sS {{\cal S}} \def\sT {{\cal T}} \def\sU {{\cal U}}
\def\sV {{\cal V}} \def\sW {{\cal W}} \def\sX {{\cal X}}
\def\sY {{\cal Y}} \def\sZ {{\cal Z}}

\def\bA {{\mathbb A}} \def\bB {{\mathbb B}} \def\bC {{\mathbb C}}
\def\bD {{\mathbb D}} \def\bE {{\mathbb E}} \def\bF {{\mathbb F}}
\def\bG {{\mathbb G}} \def\bH {{\mathbb H}} \def\bI {{\mathbb I}}
\def\bJ {{\mathbb J}} \def\bK {{\mathbb K}} \def\bL {{\mathbb L}}
\def\bM {{\mathbb M}} \def\bN {{\mathbb N}} \def\bO {{\mathbb O}}
\def\bP {{\mathbb P}} \def\bQ {{\mathbb Q}} \def\bR {{\mathbb R}}
\def\bS {{\mathbb S}} \def\bT {{\mathbb T}} \def\bU {{\mathbb U}}
\def\bV {{\mathbb V}} \def\bW {{\mathbb W}} \def\bX {{\mathbb X}}
\def\bY {{\mathbb Y}} \def\bZ {{\mathbb Z}}
\def\R {{\mathbb R}} \def\RR {{\mathbb R}} \def\H {{\mathbb H}}
\def\n{{\bf n}} \def\Z {{\mathbb Z}}

\newcommand{\expr}[1]{\left( #1 \right)}
\newcommand{\cl}[1]{\overline{#1}}
\newtheorem{thm}{Theorem}[section]
\newtheorem{lemma}[thm]{Lemma}
\newtheorem{defn}[thm]{Definition}
\newtheorem{prop}[thm]{Proposition}
\newtheorem{corollary}[thm]{Corollary}
\newtheorem{remark}[thm]{Remark}
\newtheorem{example}[thm]{Example}
\numberwithin{equation}{section}
\def\ee{\varepsilon}
\def\qed{{\hfill $\Box$ \bigskip}}
\def\NN{{\mathcal N}}
\def\AA{{\mathcal A}}
\def\MM{{\mathcal M}}
\def\BB{{\mathcal B}}
\def\CC{{\mathcal C}}
\def\LL{{\mathcal L}}
\def\DD{{\mathcal D}}
\def\bD{{\mathbb D}}
\def\FF{{\mathcal F}}
\def\EE{{\mathcal E}}
\def\QQ{{\mathcal Q}}
\def\SS{{\mathcal S}}
\def\RR{{\mathbb R}}
\def\R{{\mathbb R}}
\def\L{{\bf L}}
\def\K{{\bf K}}
\def\S{{\bf S}}
\def\A{{\bf A}}
\def\E{{\mathbb E}}
\def\bG{{\mathbb G}}
\def\F{{\bf F}}
\def\P{{\mathbb P}}
\def\N{{\mathbb N}}
\def\eps{\varepsilon}
\def\wh{\widehat}
\def\wt{\widetilde}
\def\pf{\noindent{\bf Proof.} }
\def\pff{\noindent{\bf Proof} }
\def\cp{\mathrm{Cap}}

\newcommand{\dc}{\mathfrak{A}}
\newcommand{\X}{\mathfrak{X}}
\newcommand{\dom}{\mathcal{D}}
\newcommand{\pr}{\mathbf{P}}
\newcommand{\ex}{\mathbf{E}}
\newcommand{\M}{\mathcal{M}}
\newcommand{\sub}{\subseteq}
\newcommand{\ph}{\varphi}
\newcommand{\ro}{\varrho}
\newcommand{\pot}{\mathcal{U}}
\newcommand{\potY}{\mathcal{V}}
\newcommand{\norm}[1]{\left\| #1 \right\|}
\newcommand{\snorm}[1]{\| #1 \|}
\newcommand{\set}[1]{\left\{ #1 \right\}}
\newcommand{\abs}[1]{\left| #1 \right|}
\newcommand{\as}{\text{-a.s.}}
\newcommand{\aevery}{\text{-a.e.}}

\newcommand{\itref}[1]{\ref{#1}}

\title{\Large \bf
Scale invariant boundary Harnack principle at infinity for Feller processes
}

\author{{\bf Panki Kim}\thanks{This work was  supported by the National Research Foundation of
Korea(NRF) grant funded by the Korea government(MSIP) (No. NRF-2015R1A4A1041675)
}
\quad {\bf Renming Song\thanks{Research supported in part by a grant from
the Simons Foundation (208236)}} \quad and
\quad {\bf Zoran Vondra\v{c}ek}
\thanks{Research supported in part by the Croatian Science Foundation under the project 3526}
}

\date{}

\maketitle

\begin{abstract}
In this paper we prove a uniform and scale invariant
boundary Harnack principle at infinity for a large class
of purely discontinuous Feller processes in metric measure spaces.
\end{abstract}

\noindent {\bf AMS 2010 Mathematics Subject Classification}: Primary 60J50, 31C40; Secondary 31C35, 60J45, 60J75.

\noindent {\bf Keywords and phrases:} Boundary Harnack principle, harmonic functions,
purely discontinuous Feller process

\section{Introduction}\label{s:intro}
The boundary Harnack principle (BHP) is a result roughly saying that non-negative functions, which are harmonic in an open
set and vanish near a portion of the boundary of that open set,
have the same boundary decay rate near that portion of the boundary. The BHP
was first proved independently in \cite{anc, dah, wu} for classical harmonic functions in Lipschitz domains.
Since then, it has been extended to more general
diffusions and more general domains.

In \cite{Bogdan97}, the BHP was established for harmonic functions
of symmetric $\alpha$-stable processes, $\alpha\in (0, 2)$, in Lipschitz domains.
This was the first BHP for discontinuous Markov processes. Since then,
the result of \cite{Bogdan97} has been generalized in various directions. \cite{SW99}
extended it to harmonic functions of symmetric $\alpha$-stable
processes in $\kappa$-fat open sets, with the constant depending on the local geometry
near the boundary. A uniform version of it was established in \cite{BKK08} for harmonic
functions of symmetric $\alpha$-stable processes in arbitrary open sets. The BHP
of \cite{BKK08} is uniform in the sense that the constant does not
depend on the open set itself.
Note that such uniform version does not hold for Brownian motion.

In another direction, the BHP has been
generalized to different classes of discontinuous Markov processes. For example,
it was extended to a large class
of subordinate Brownian motions in \cite{KSV09, KSV12yan}.
In \cite{KSV12} the uniform BHP was extended to a large class of rotationally symmetric
L\'evy processes and in \cite{KM2} it was extended to a class of subordinate
Brownian motions including geometric stable processes. The main result of \cite{KSV12}
has been extended to a large class of symmetric L\'evy processes in \cite{KSV15}.
A BHP with explicit decay rate was established in \cite{KSV12plms, KSV14a} for a large class of subordinate
Brownian motions in $C^{1, 1}$ open sets.
For BHP with respect to subordinate Brownian motions with Gaussian components,
see \cite{CKSV, KSV13}.

Recently, a very general BHP for discontinuous Feller processes in metric measure spaces 
has  been proved in \cite{BKuK} under some comparability assumptions on
the jump kernel and a Urysohn-type property of the domain of the generator of the process.
The main result of \cite{BKuK} is not scale invariant in general. It was shown in
\cite{BKuK} that, under a stable-like scaling condition, a scale invariant BHP holds.

All the BHPs mentioned above deal with the decay of harmonic functions
near finite boundary points. In the case of symmetric $\alpha$-stable
processes, by using the inversion with respect to spheres, the Kelvin transform and
the BHP near finite boundary points,
\cite{Kwa} obtained a BHP at infinity for harmonic functions
in unbounded open sets.  The argument using inversion with respect to spheres and
the Kelvin transform does not work for more general L\'evy processes.
By using a different, more involved argument, a BHP at infinity was established in \cite{KSV14}
for a large class of symmetric L\'evy processes under a global weak scaling condition
on the L\'evy exponents.

Motivated by the result and the method from \cite{BKuK},
in this paper we prove a uniform and scale invariant BHP at infinity for a class of
purely discontinuous Feller processes in metric measure spaces.
Even in the special case of symmetric L\'evy processes, the BHP at infinity of this
paper is more general than that of \cite{KSV14}
since we will only assume that the L\'evy exponents satisfy a
weak scaling condition near the origin.
We will also give a uniform and scale invariant BHP near finite boundary points.

We start the paper by recalling  the setting and basic assumptions of \cite{BKuK}. Let $(\X, d)$ be a metric space such that  all bounded closed sets are compact and let $m$ be a $\sigma$-finite measure on $\X$ with full support. Let $R_0 \in (0, \infty]$  (the localization radius of $(\X, d)$) be such that $\X \setminus B(x, 2 r) \neq \emptyset$ for all $x \in \X$ and all $r < R_0$. We will consider a large class of Feller processes $X=(X_t,t\ge0; \P_x, x\in \X; \mathcal{F}_t, t\ge0)$ on $\X$ satisfying several assumptions. The first assumption is strong duality and Hunt's hypothesis (H).

\smallskip
\noindent
{\bf A}: $X$ is a Hunt process admitting a strong dual process $\widehat{X}$ with respect to the measure $m$ and $\widehat{X}$ is also a Hunt process.
The transition semigroups $(P_t)$ and $(\widehat{P}_t)$ of $X$ and $\widehat{X}$ are both Feller and
strongly Feller.
Every semi-polar set of $X$ is polar.

\smallskip
In the sequel, all objects related to the dual process $\widehat{X}$ will be denoted by a hat.
Recall that a set is polar (semi-polar, respectively) for $X$ if and only if it is polar (semi-polar, respectively) for
$\wh X$ (see \cite[VI. (1.19)]{BG}).
Under assumption {\bf A} the process $X$ admits a (possibly infinite) Green function $G(x,y)$ serving as a density of the occupation measure: $G(x,A):=\E_x \int_0^{\infty}\ind_{(X_t\in A)} dt =\int_A G(x,y)m(dy)$. Moreover, $G(x,y)=\widehat{G}(y,x)$ for all $x,y\in \X$,
cf.~\cite[VI.1]{BG}.
Further, if $D$ is an open subset of $\X$ and $\tau_D=\inf\{t>0:\, X_t\notin D\}$ the exit time from $D$, the killed process $X^D$ is defined by $X_t^D=X_t$ if $t<\tau_D$ and $X_t^D=\partial$
where $\partial$ is an extra point added to $\X$.
Then $X^D$ admits a unique (possibly infinite) Green function (potential kernel) $G_D(x,y)$ such that for every non-negative Borel function $f$,
$$
G_D f(x):=\int_{D} f(y)G_D(x,y)dy=\E_x \int_0^{\tau_D}f(X_t)dt\,  ,
$$
and $G_D(x,y)=\widehat{G}_D(y,x)$, $x,y\in D$, with $\widehat{G}_D(y,x)$ the Green function of $\widehat{X}^D$.
For the details we refer the readers to \cite[pp.480--481]{BKuK} and the references therein.
We say $D$ is Greenian if the Green function $G_D(x,y)$ is finite for all $x,y\in D$, $x\neq y$. Under this assumption the process $X^D$ is transient in the sense that there exists a non-negative Borel function $f$ on $D$ such that $0<G_D f<\infty$ (and the same is true for $\widehat{X}$).

Let $C_0(\X)$ stand for the Banach space of bounded continuous functions on $\X$ vanishing at infinity.
Let $\sA$ and $\widehat \sA$ be the generators of $(P_t)$ and $(\widehat{P}_t)$ in $C_0(\X)$ respectively. The second assumption is a Urysohn-type condition.

\smallskip
\noindent
{\bf B}: There is a linear subspace $\dom$ of $\dom(\sA) \cap \dom(\widehat{\sA})$ satisfying the following condition:
For any compact $K$ and open $D$ with $K\subset D\subset \X$, the collection $\dom(K, D)$ of functions $f\in \dom$ satisfying the conditions (i) $f(x)=1$ for $x\in K$; (ii) $f(x)=0$ for $x\in \X\setminus D$; (iii) $0\le f(x)\le 1$ for $x\in \X$, and (iv) the boundary of the set $\{x \, : \, f(x) > 0\}$ has
zero $m$ measure, is nonempty. We let
$$
\rho(K, D):=\inf_{f\in \dom(K, D)} \sup_{x \in \X} \max(\sA f(x), \widehat{\sA} f(x)).
$$
\smallskip

Assumption {\bf B} implies that the jumps of $X$ satisfy the following L\'evy system formula:
for every stopping time $T$,
\begin{equation}\label{e:lsformula}
\E_x\sum_{s\in (0, T]}f(X_{s-}, X_s)= \E_x\int^T_0\int_{\X}f(X_s, z)J(X_s, dz)ds.
\end{equation}
Here $f:\X\times\X\to[0, \infty]$, $f(x, x)=0$ for all $x\in \X$,
and $J$ is a kernel on $\X$ (satisfying $J(x, \{x\})=0$ for all $x\in\X$), called the L\'evy kernel of $X$. As a consequence, the following
Ikeda-Watanabe type formula is valid:
\begin{equation}\label{e:ikeda-watanabe}
\P_x(X_{\tau_D}\in E, X_{\tau_D-}\neq X_{\tau_D}, \tau_D<\zeta)=\int_D G_D(x,dy)J(y,E), \quad x\in D, E\subset \X\setminus D\,
\end{equation}
where $\zeta$ is the life time of $X$.
Furthermore, the L\'evy kernel $J$ satisfies
\begin{equation}\label{e:levykernel}
Jf(x):=\int_{\X}f(y)J(x,dy)=\lim_{t\downarrow0}\frac{\E_xf(X_t)}{t}
\end{equation}
for all bounded continuous function $f$ on $\X$ and $x\in \X\setminus\mbox{supp} f$.
The L\'evy kernel $\widehat{J}(y, dx)$ of $\widehat{X}$ is defined in a similar manner. By duality, $J(x, dy)m(dx)=\widehat{J}(y, dx)m(dy)$. Further, it follows from \eqref{e:levykernel} that if $f\in \dom(\sA)$ and $x\in \X\setminus \mathrm{supp}(f)$, then $Jf(x)=\sA f(x)$.
Again, for these facts we refer the reader to \cite[p.482]{BKuK} and the reference therein.

Our next assumption is only a part of the corresponding assumption in \cite{BKuK}.

\smallskip
\noindent
{\bf C}: The L\'evy kernels of $X$ and $\widehat{X}$ have the form $j(x, y)m(dy)$ and $\widehat{j}(x, y)m(dy)$ respectively, where $j(x, y)=\widehat{j}(y, x)>0$
for all $x, y\in\X, x\neq y$.

\smallskip
For an open set $D\subset \X$, let
\begin{equation}\label{e:poisson-kernel}
P_D(x,z):=\int_D G_D(x,y)j(y,z)m(dy), \qquad x\in D, z\in D^c,
\end{equation}
be the Poisson kernel of $D$ with respect to $X$. It follows from \eqref{e:ikeda-watanabe} and \eqref{e:poisson-kernel} that $P_D(x,z)$ is the exit density of $X$
from $D$ through jumps:
$$
\P_x(X_{\tau_D}\in E, X_{\tau_D-}\neq X_{\tau_D}, \tau_D<
\zeta)=\int_E P_D(x,z)m(dz), \quad x\in D, E\subset \X\setminus D\,  .
$$

Assumptions {\bf A}, {\bf B}  and {\bf C} will be in force throughout this  paper.
In the next section we will assume that the localization radius $R_0=\infty$, that $X$ and $\widehat{X}$ are conservative, and will add assumptions needed in order to study the behavior of non-negative harmonic functions at infinity. Our main result is a scale invariant approximate factorization of
non-negative function harmonic in unbounded open set,
Theorem \ref{t:main-infty},
from which the scale invariant uniform boundary Harnack principle, Corollary \ref{c:bhp-infty}, immediately follows.
Proofs of these results will be given in Section \ref{s:proofs}. In Section \ref{s:finite},
we introduce a different additional set of assumptions and state
the scale invariant uniform boundary Harnack principle at a finite boundary point.
Proofs of the results in Section \ref{s:finite} are deferred to the Appendix.
In Section 5 we discuss examples of processes which satisfy the assumptions of this paper. These include some symmetric and isotropic L\'evy processes, strictly stable (not necessarily symmetric) processes in $\R^d$, processes obtained by subordinating a Feller diffusion on unbounded Ahlfors regular $n$-sets, and some space non-homogeneous processes on $\R^d$.
Finally, in Section \ref{s:green} we study boundary behavior of the Green function $G_D$ at the regular boundary points.
These results will be used in subsequent papers.

Notation: We will use the following conventions in this paper.
 $c, c_0,
c_1, c_2, \cdots$ stand for constants
whose values are unimportant and which may change from one
appearance to another. All constants are positive finite numbers.
The labeling of the constants $c_0, c_1, c_2, \cdots$ starts anew in
the statement of each result. We will use ``$:=$"
to denote a definition, which is  read as ``is defined to be".
We denote $a \wedge b := \min \{ a, b\}$.
The notation $f\asymp g$  means that the quotient $f(t)/g(t)$ stays bounded
between two positive numbers on their common domain of definition.  For $x\in \X$ and $r>0$ we denote by $B(x,r)$ be the open ball centered at $x$ with radius $r$ and by $\overline{B}(x,r)$ the closure of $B(x,r)$. Further, for $0<r<R$, let $A(x,r,R)=\{y\in \X: r <d(x,y)<R\}$ be the open annulus around $x$, and $\overline{A}(x,r,R)$ the closure of $A(x,r,R)$.
Throughout the paper we will adopt the convention that
$X_{\zeta}=\partial$ and $u(\partial)=0$ for every function $u$.

\section{Additional assumptions and main result}\label{s:main-result}
We recall that assumptions {\bf A}, {\bf B}  and {\bf C}  are in force throughout the paper.
In order for the boundary Harnack principle at infinity to make sense the ambient space $\X$ must be unbounded. Hence in this and the next section we assume that the localization radius $R_0=\infty$. We further assume that both $X$ and $\widehat{X}$ are conservative processes: For every $t\ge 0$, $P_t 1=\widehat{P}_t 1=1$.

From now on we fix a point $z_0$ in $\X$ which will serve as the center of the space. For any $r>0$, let
$$
V(r)=V(z_0,r):=m(B(z_0,r))
$$
denote the volume of the ball of radius $r$ centered at $z_0$.
We assume that $V:[0,\infty)\to [0,\infty)$ satisfies the following two properties: 
(i) The doubling property:  There exists
$c>1$ such that
\begin{equation}\label{e:V-doubling}
V(2r)\le c V(r), \quad r>0\,  ,
\end{equation}
and (ii) There exist $c>1$, $r_0>0$ and $n_0 \in \N$ with $n_0 \ge 2$ such that
\begin{equation}\label{e:V-another}
V(n_0r) \ge c V(r), \quad r \ge r_0\,  .
\end{equation}

We further assume the existence of a
non-decreasing
 function $\Phi=\Phi(z_0, \cdot):[0,\infty)\to [0,\infty)$ which satisfies the doubling property: There exists
 $c>1$ such that
\begin{equation}\label{e:Phi-doubling}
\Phi(2r)\le c \, \Phi(r), \quad r>0\,  .
\end{equation}
The function $\Phi$ will play a crucial role in obtaining scale invariant results. Examples of such functions will be given in Section \ref{s:examples}. At the moment it suffices to say that
in case of isotropic L\'evy process in $\R^d$,
we have that $\Phi(r)=1/\Psi(r^{-1})$,
where $x \mapsto \Psi(|x|)$ is the L\'evy exponent of the process.

Let $C_\infty(\X)$ be the Banach space of continuous functions $f$ on $\X$ such that $f$ has a limit at infinity. We will use $\|\cdot\|$ to denote the sup norm. It is obvious that any function $f\in C_\infty(\X)$ is the sum of function in $C_0(\X)$ and a constant. It is well known that the semigroup of $X$ being Feller is equivalent to the following conditions: (i) for any $f\in C_\infty(\X)$, $P_t f\in C_\infty(\X)$; (ii) for any $f\in C_\infty(\X)$, $\lim_{t\to 0}\|P_tf-f\|=0$. We will also use $\sA$ (respectively $\widehat{\sA}$) to denote the generator of $(P_t)$ (respectively $(\widehat{P}_t)$) in $C_\infty(\X)$. It follows easily from the conservativeness of $X$
that constant functions are in $\dom(\sA)$ and that, for any constant $c$, $\sA c=0$.

We are now ready for some additions to assumptions {\bf B} and {\bf C} and an additional assumption.
In the following assumptions, $r_0$ is a positive number.
Recall the notation $\dom(K,D)$ from assumption {\bf B}.

\smallskip
\noindent
{\bf B2-a}$(z_0, r_0)$:
For any $a\in (1,2]$, there exists $c=c(z_0,a)$ such that for any $r\ge r_0$,
\begin{align*}
 \ro(r) & :=  \inf_{f\in \dom(\overline{B}(z_0,r ), B(z_0, ar)) } \sup_{x \in \X} \max(\sA(1-f)(x), \widehat{\sA}(1- f)(x))\\
&= \inf_{f\in \dom(\overline{B}(z_0, r), B(z_0, ar)) } \sup_{x \in \X} \max(-\sA f(x), -\widehat{\sA} f(x))\le  \frac{c}{\Phi(r)}.
\end{align*}

\noindent
{\bf B2-b}$(z_0, r_0)$:
For any $a\in (1,2]$, there exists $c=c(z_0,a)$ such that for any $r\ge r_0$ and any $f\in \dom(\overline{B}(z_0, r), B(z_0, ar))$,
\begin{align*}
 \max(\sA f(x), \widehat{\sA} f(x)) \le  c V(r) j(x,z_0),  \quad x \in \overline{A}(z_0, r, (a+1)r).
\end{align*}

To assumption {\bf C} we add

\smallskip
\noindent
{\bf C2}$(z_0, r_0)$:
For any $a\in (1,2]$,  there exists $c= c(z_0, a)$
such that for $r\ge r_0$, $x \in B(z_0, r)$ and $y \in \X \setminus B(z_0, ar)$,
\begin{align}\label{e:AC2}
 c^{-1}{j(z_0, y)} \le {j(x, y)} \le c{j(z_0, y)} , &&
 c^{-1}\, {\widehat j(z_0, y)} \le {\widehat j(x, y)} \le c{\widehat j(z_0, y)},
\end{align}
and
\begin{equation}\label{e:AC21}
\inf_{y\in \overline{A}(z_0,r,ar)} \min(j(z_0,y), \widehat{j}(z_0,y))\ge \frac{c}{V(r)\Phi(r)}\, .
\end{equation}
\smallskip

Note that by assumptions {\bf B2-b}$(z_0, r_0)$ and {\bf C2}$(z_0, r_0)$, and the doubling property of $V$,
any $f\in \dom(\overline{B}(z_0, r), B(z_0, ar))$ satisfies
\begin{align}\label{e:Assumptions B2(b)}
 \max(\sA f(x), \widehat{\sA} f(x)) \le c  {\bf 1}_{B(z_0, r)^c} (x)  V(r) j(x,z_0),
 \qquad r\ge r_0,
\end{align}
for a constant $c=c(z_0,a)>0$.
In fact,  for $f \in  \dom(\sA) \cap \dom(\widehat{\sA})$ such that $f(x) = 1$ for $x \in   B(z_0, r)$, $f(x) = 0$ for $x \in \X \setminus \overline{B}(z_0, ar)$ and $0 \le f(x) \le 1$ for $x \in \X$,
we have, using  the doubling property of $V$,
\begin{align}
&\sA f(x) {\bf 1}_{\X \setminus \overline{B}(z_0, (a+1)r)} (x)  =
{\bf 1}_{\X \setminus \overline{B}(z_0, (a+1)r)} (x)\int_{B(z_0, ar)} f(y) j(x,y)m(dy) \nonumber\\
&\le c {\bf 1}_{\X \setminus \overline{B}(z_0, (a+1)r)} (x) V(ar)j(x,z_0)\le c_2{\bf 1}_{\X \setminus \overline{B}(z_0, (a+1)r)} (x) V(r)j(x,z_0) .
\end{align}
Combining this with assumption {\bf B2-b}$(z_0, r_0)$, we get \eqref{e:Assumptions B2(b)}.

Our final assumption concerns Green functions of complements of balls.

\smallskip
\noindent
{\bf D2}$(z_0, r_0)$:
$\overline{B}(z_0, r_0)^c$ is Greenian.
For every $a\in (1,2)$, there exists a constant $c=c(z_0,a)$ such that for all $r\ge r_0$,
 $$
\sup_{x \in \overline{B}(z_0, ar)} \sup_{y \in  \overline{A}(z_0, 2r, 4r)} \max(G_{\overline{B}(z_0, r)^c}(x, y), \widehat{G}_{\overline{B}(z_0, r)^c}(x, y)) \le c\frac{\Phi(r)}{V(r)}.
$$

Recall that a non-negative function $u:\X\to [0,\infty)$ is said to be regular harmonic in an open set $D\subset \X$ if
$$
u(x)=\E_x \left(u(X_{\tau_D})\right), \quad \textrm{for all }x\in D\,  .
$$
By the strong Markov property, the equality above holds for every stopping time $\tau\le \tau_D$.

Recall also that for an open set $D\subset \X$, a point $x\in \partial D$ is said to be regular for $D^c$
with respect to $X$ if $\P_x(\tau_D=0)=1$. Let $D^{\mathrm{reg}}$ denote the set of points $x\in \partial D$
which are regular for $D^c$ with respect to $X$.

Now we can state our main theorem.
\begin{thm}\label{t:main-infty}
Assume $z_0\in\X$.
Suppose that, in addition to {\bf A}, {\bf B} and {\bf C},
assumptions \eqref{e:V-doubling}--\eqref{e:Phi-doubling},
{\bf B2-a}$(z_0, r_0)$,
{\bf B2-b}$(z_0, r_0)$, {\bf C2}$(z_0, r_0)$ and {\bf D2}$(z_0, r_0)$  hold true for some $r_0>0$.
For any $a\in (1,2)$, there exists $C_1=C_1(z_0, a)>1$ such that for any $r \ge r_0$, any open set
$D \subset \overline{B}(z_0,r)^c$ and any non-negative function $u$ on $\X$
which is regular harmonic with respect to $X$ in $D$ and vanishes on $\overline{B}(z_0,r)^c\cap\left(\overline{D}^c\cup D^{\mathrm{reg}}\right)$,
it holds that
\begin{equation}\label{e:bhp-inf}
C_1^{-1}P_{D}(x,z_0) \int_{B(z_0,2ar)}u(z)\, m(dz) \le u(x) \le C_1 P_{D}(x,z_0) \int_{B(z_0, 2ar)}u(z)\, m(dz)
\end{equation}
for all $x\in D\cap \overline{B}(z_0,8r)^c$.
\end{thm}

As a consequence of Theorem \ref{t:main-infty}, one immediately gets the following scale
invariant uniform boundary Harnack principle at infinity.

\begin{corollary}[Boundary Harnack Principle at Infinity]\label{c:bhp-infty}
Let $z_0\in
\X$. Assume that, in addition to {\bf A}, {\bf B} and {\bf C}, assumptions
\eqref{e:V-doubling}--\eqref{e:Phi-doubling},
{\bf B2-a}$(z_0, r_0)$,
{\bf B2-b}$(z_0, r_0)$, {\bf C2}$(z_0, r_0)$ and {\bf D2}$(z_0, r_0)$  hold true for some $r_0>0$.
There exists $C_2=C_2(z_0)>1$ such that for any $r \ge 2r_0$, any open
set $D \subset \overline{B}(z_0,r)^c$ and any non-negative
functions $u$ and $v$ on $\X$ which are regular harmonic with respect to $X$ in $D$ and
vanish on $\overline{B}(z_0,r)^c\cap\left(\overline{D}^c\cup D^{\mathrm{reg}}\right)$, it holds that
\begin{equation}\label{e:bhp-inf-cor}
    C_2^{-1} \frac{u(y)}{v(y)} \le  \frac{u(x)}{v(x)}\le C_2 \frac{u(y)}{v(y)}\, ,\qquad \textrm{for all }x,y\in D\cap \overline{B}(z_0,8r)^c\,  .
\end{equation}
\end{corollary}

\begin{remark}\label{r:symmetry}{\rm
Note that all our assumptions are symmetric in $X$ and $\widehat{X}$. Therefore, Theorem \ref{t:main-infty} and Corollary \ref{c:bhp-infty} hold for co-harmonic functions as well.
}
\end{remark}

\section{Proofs}\label{s:proofs}
Throughout this section, $z_0$ is a fixed point in $\X$.
We will always  assume in this section that assumptions  {\bf A}, {\bf B}, {\bf C}, \eqref{e:V-doubling}--\eqref{e:Phi-doubling},
{\bf B2-a}$(z_0, r_0)$, {\bf B2-b}$(z_0, r_0)$, {\bf C2}$(z_0, r_0)$ and {\bf D2}$(z_0, r_0)$  hold true
for some $r_0>0$ and give a proof of Theorem \ref{t:main-infty}.

\begin{prop}\label{p:pcom}
There exists a constant $c>0$ such that for every  $r \ge r_0$
and $D \subset \X \setminus \overline{B}(z_0, r)$,
$$
P_D(x, z_0) \le c\, \frac{1}{V(r)}, \quad x \in D.
$$
\end{prop}
\pf
By \eqref{e:AC2} and the fact that $j(x,y)=\widehat j(y,x)$, for every $y \in {B}(z_0, r/2)$ and $x \in D$, we have
\begin{align*}
&P_D(x, z_0) \le P_{\overline{B}(z_0, r)^c}(x, z_0)
=\int_{\overline{B}(z_0, r)^c}G_{\overline{B}(z_0, r)^c}(x,z) j(z, z_0)m(dz)\\
&\le c_1 \int_{\overline{B}(z_0, r)^c}G_{\overline{B}(z_0, r)^c}(x,z) j(z, y)m(dz)= c_1 P_{\overline{B}(z_0, r)^c}(x, y).
\end{align*}
Thus, by integrating over the ball $B(z_0,r/2)$ and using the doubling property of $V(r)$,
\begin{align*}
P_D(x, z_0) \le \frac{c_2}{V(r/2)} \int_{{B}(z_0, r/2)} P_{\overline{B}(z_0, r)^c}(x, y) m(dy) \\
\le \frac{c_3}{V(r)} \int_{{B}(z_0, r/2)} P_{\overline{B}(z_0, r)^c}(x, y) m(dy)
\le  \frac{c_3}{V(r)} .
\end{align*}
\qed

Let $a\in (1,2)$. For each $r\ge r_0$, we consider a function $\ph^{(r)}\in \dom(\overline{B}(z_0, r), B(z_0, ar))$,
and let $\phi^{(r)}=1-\ph^{(r)}$ and $V^{(r)}=\{x\in \X: \phi^{(r)}(x)>0\}=\{x\in \X: \ph^{(r)}(x)<1\}$.
Note that, by choosing $\ph^{(r)}$ appropriately,  we can achieve that
$$
\delta^{(r)}:= \sup_{x \in B(z_0, ar)} \max(\sA \phi^{(r)}(x), \widehat{\sA} \phi^{(r)}(x)) \le \frac{2c}{\Phi(r)},
$$
where $c=c(z_0,a)$ is the constant in assumption {\bf B2-a}$(z_0, r_0)$.

In what follows, our analysis and results are valid for all $r \ge r_0$ with constants depending on $a\in (1,2)$,
but \emph{not} on $r$.
To ease the notation in the remaining part of the section we drop the superscript $r$ from $\ph^{(r)}$, $\phi^{(r)}$ and $V^{(r)}$ and write simply $\ph$, $\phi$ and $V$.

Let
\begin{align}\label{e:defpsi}
 \psi(x)  = \frac{\max(\sA \phi(x), \widehat{\sA} \phi(x), \delta (1 - \phi(x)))}{\phi(x)} \, ,
\quad   x \in \X ,
\end{align}
with the convention $1/0=\infty$. Note that $\psi(x)=\infty$ for $x\in V^c$, and $\psi(x)=0$ for $x\in \overline{B}(z_0,ar)^c$.
We define two right-continuous additive functionals by
\begin{align}\label{e:eAt}
A_t = \lim_{\eps \searrow 0} \int_0^{t + \eps} \psi(X_s) ds\quad \text{and} \quad
 \widehat A_t = \lim_{\eps \searrow 0} \int_0^{t + \eps} \psi( \widehat X_s) ds.
\end{align}

We follow the idea in  \cite{BKuK} to mollify the distribution of $X(\tau_{\overline{B}(z_0, r)^c})$ by letting the particle lose mass gradually,
with intensity $\psi(X_t)$, before time $\tau_{\overline{B}(z_0, r)^c}$.

We define two right-continuous strong Markov multiplicative functionals
$M_t = \exp(-A_t)$ and $\widehat M_t = \exp(-\widehat A_t)$.
By the argument on \cite[p.~492]{BKuK} and the references therein, $M$ and
$\widehat M$ are in fact exact strong Markov multiplicative functionals.
As on \cite[p.~492]{BKuK}, we consider the semigroup of operators $T^\psi_t f(x) = \E_x(f(X_t) M_t)$
associated with the multiplicative functional $M$
and the semigroup of operators $\widehat T^\psi_t f(x) = \E_x(f(\widehat X_t) \hat M_t)$
associated with the multiplicative functional $\widehat M$.
 $T^\psi_t$ is the transition operator of the subprocess $X^\psi$ of $X$.
The subprocess $\widehat{X}^\psi$ of $\widehat{X}$ corresponding to the multiplicative
functional $\widehat{M}$ is the dual process of $X^\psi$.
Thus the  potential densities of $X^\psi$ and $\widehat{X}^\psi$
satisfy $\widehat{G}^\psi(x, y) = G^\psi(y, x)$
and
\begin{align}\label{e:Gmon}
 G^\psi(x,y) \le  G_V(x,y)  \le  G_{\overline{B}(z_0, r)^c} (x,y),  \quad (x,y) \in V \times V.
\end{align}

Let $\tau_a = \inf \set{t \ge 0 : A_t \ge a}$ and
\begin{equation}\label{e:psi-tau-a}
 \pi^\psi f(x) =-\E_x \int_{[0, \infty)} f(X_t) d M_t = \E_x \expr{\int_0^\infty f(X_{\tau_a}) e^{-a} da} = \int_0^\infty \E_x(f(X_{\tau_a})) e^{-a} da,
\end{equation}
where the second equality follows by substitution.
By following the arguments in \cite[pp.~492--493]{BKuK} line by line,
one can see that  $\pi^\psi f$ can be written in the following two ways:
If $f$ is nonnegative and vanishes in $\X \setminus (B(z_0, ar) \cap V)$, then
\begin{align}\label{e:pipsi1}
 \pi^\psi f(x) & = G^\psi (\psi f)(x) = \int_{V \cap B(z_0, ar)} G^\psi(x, y)
 \psi(y) f(y) m(dy), \quad
 x \in V,
\end{align}
and if $f \in \dom(\sA)$ vanishes in $V$, then
for all $x \in V$,
\begin{align}\label{e:pipsi2}
& \pi^\psi f(x)  = G^\psi \sA f (x) = \int_V G^\psi(x, y) \sA f(y) m(dy)\nonumber\\
&=\int_V G^\psi(x, y) \int_{\X \setminus V} f(z) j(y, z)  m(dz) m(dy)\nonumber\\
 &=\int_{\X \setminus V} \expr{\int_V G^\psi(x, y) j(y, z) m(dy)} f(z) m(dz).
\end{align}

Since Lemmas 4.6--4.7 and Corollary 4.8 of \cite{BKuK} and their proofs  hold under our setting, we have
\begin{align} \label{e:Vzero}
\pi^\psi(x, \partial V) = 0, \qquad x\in V.
\end{align}

Repeating the argument of the proof of \cite[Lemma 4.10]{BKuK}, we get that if $f$ is regular harmonic
in $\overline{B}(z_0, r)^c$ with respect to $X$,
then  $f(x) = \pi^\psi f(x)$
for all $x\in B(z_0, 2r)^c$.
The main step of the proof is to get the correct estimate of $\pi^\psi(x, dy)/m(dy)$.

We recall the following notation and results from \cite{BKuK}.
Let $U$ be an open subset of $V$. For any nonnegative or bounded $f$ and $x \in V$, we let
\begin{align*}
 \pi^\psi_U f(x) & = \E_x (f(X_{\tau_U}) M_{\tau_U-}) , & G^\psi_U f(x) & = \E_x \int_0^{\tau_U} f(X_t) M_t dt .
\end{align*}
 $G^\psi_U$ admits a density $G^\psi_U(x, y)$, and we have $G^\psi_U(x, y) \le G_U(x, y)$, $G^\psi_U(x, y) \le G^\psi(x, y)$.
 For any $f \in \dom(\sA)$, we have
\begin{align}
\label{e:psid2}
 \pi^\psi_U f(x) & = G^\psi_U (\sA - \psi) f(x) + f(x) , && x \in V .
\end{align}
In particular, by an approximation argument,
\begin{align}
\label{e:psilevy}
  \pi^\psi_U(x, E) & = \int_U G^\psi_U(x, y) J(y, E) m(dy) , && x \in U , \, E \sub \X \setminus \overline{U} .
\end{align}

By the definition of $\psi$, we have that $(\sA - \psi) \phi(x) \le 0$ for $x \in \X$ (and, in particular, for all $x\in V$). Thus using \eqref{e:psid2},
the proof of the next result is the same as that of \cite[Lemma 4.4]{BKuK}.

\begin{lemma}
\label{l:4.4}
Let $U = V \cap {B}(z_0, ar)$. Then
\begin{align}
\label{e:4.4}
 \pi^\psi_U(x, \overline{B}(z_0, ar)^c) & \le \phi(x) , && x \in U .
\end{align}
\end{lemma}

Recall that $n_0$ is the natural number in \eqref{e:V-another}.
It follows from the doubling properties of $\Phi$ and $V$ that  with $c_0>1$
\begin{equation}\label{e:dbP1}
\Phi(n_0^k r)\le c_0^k \Phi(r) \quad \text{for all } k\ge 1, r>0,
\end{equation}
\begin{equation}\label{e:dbP2}
V(n_0r/2)\Phi(n_0r/2) \le c_0V(r)\Phi(r) \quad \text{for all }   r >0.
\end{equation}
Let $s\ge r_0$.
Thus applying \eqref{e:AC21}  to each annulus and the monotonicity of $\Phi$ and $V$
we have that
 for $n_0^k s < d(z, z_0)\le n_0^{k+1}s$,
\begin{equation}\label{e:dbP3}
j(z_0,z) \ge
\frac{c_1}{V(n_0^{k+1} s/2)\Phi(n_0^{k+1} s/2)} \ge \frac{c_1c_0^{-1}}{V(n_0^k s)\Phi(n_0^k s)}.
\end{equation}
In the second inequality above we have used \eqref{e:dbP2}.
   Let $s\ge r_0$. By using  \eqref{e:dbP3}
 in the second line, the assumption \eqref{e:V-another} (with the constant $c_2>1$) and \eqref{e:dbP1} in the last line, we have
\begin{eqnarray*}
\int_{\overline{B}(z_0,s)^c} j(z_0,z)m(dz)&=&\sum_{k=0}^{\infty}\int_{n_0^k s < d(z, z_0)\le n_0^{k+1}s} j(z_0,z)m(dz)\\
&\ge & \sum_{k=0}^{\infty}\int_{n_0^k s < d(z, z_0)\le n_0^{k+1}s} \frac{c_1c_0^{-1}}{V(n_0^k s)\Phi(n_0^k s)} m(dz)\\
&=& c_1c_0^{-1} \sum_{k=0}^{\infty} \frac{V(n_0^{k+1}s)-V(n_0^k s)}{V(n_0^k s)\Phi(n_0^k s)}\\
&=& c_1 c_0^{-1}\sum_{k=0}^{\infty} \left(\frac{V(n_0^{k+1}s)}{V(n_0^k s)}-1\right) \frac{1}{\Phi(n_0^k s)}\\
&\ge & c_1 (c_2-1) \sum_{k=0}^{\infty} c_0^{-k-1}\frac{1}{\Phi(s)}=: c_3 \frac{1}{\Phi(s)},
\end{eqnarray*}
with $c_2>1$ and $c_1, c_3>0$ independent of $s$. In particular,  again by using \eqref{e:Phi-doubling},
\begin{equation}\label{e:jlower}
\int_{\overline{B}(z_0, br)^c} j(z_0, z)m(dz)   \ge c _4 \frac{1}{\Phi(r)}, \quad \textrm{for every }b\in (1,2] \textrm{ and all }r\ge r_0\,  ,
\end{equation}
with $c_4>0$ independent of $r$ and $b$.

\begin{lemma}
\label{l:etauU}
Let $b\in (a,2)$ and set $U = V \cap {B}(z_0, ar)$.  There exists a constant $c=c(z_0,b/a)>0$ such that
\begin{align}
\label{e:etauU}
 \int_U G^\psi_U(x, y) m(dy)  & \le c
 \phi(x) \Phi(r),
  \qquad  x \in U .
\end{align}
\end{lemma}
\pf
By Lemma \ref{l:4.4},
$
\phi(x) \ge  \pi^\psi_U(x, \overline{B}(z_0, ar)^c) \ge \pi^\psi_U(x, \overline{B}(z_0, br)^c)\,  .
$
Thus, using \eqref{e:psilevy},
$$
\phi(x) \ge \int_{\overline{B}(z_0, br)^c}  \int_U G^\psi_U(x, y) j(y, z) m(dy)m(dz)\,  .
$$
Note that, by \eqref{e:AC2}, we have
$
j(y, z) \ge c_1 j(z_0, z)$ for all  $(y,z) \in B(z_0, ar)  \times \overline{B}(z_0, br)^c\
$
with $c_1=c_1(z_0,b/a)$. Therefore, using \eqref{e:jlower} we conclude that
$$
\phi(x) \ge c_1 \int_{\overline{B}(z_0, br)^c}  j(z_0, z) m(dz) \int_U G^\psi_U(x, y)m(dy) \ge \frac{c_2}{\Phi(r)} \int_U G^\psi_U(x, y)m(dy)\, .
$$
\qed

The following Lemma is analogous  to \cite[Lemma 4.5]{BKuK}.
\begin{lemma}\label{l:4.5}
Let $b\in (a,2)$. There exists a constant  $c=c(a,b)>0$ such that
\begin{align*}
 G^\psi(x, y) & \le c  \frac{\Phi(r)}{V(r)} \phi(x), \quad   x \in V \cap \overline{B}(z_0, br) , \, y \in \overline{A}(z_0, 2r, 4r) \, .
\end{align*}
\end{lemma}
\pf
If $x \in \overline{A}(z_0, ar, br)$, then $\phi(x)=1$. Thus,
by \eqref{e:Gmon} and assumption {\bf D2}$(z_0, r_0)$,
for $x \in \overline{A}(z_0, ar, br)\subset \overline{B}(z_0, br)$ and $y \in \overline{A}(z_0, 2r, 4r)$,
\begin{align*}
 G^\psi(x, y)  \le  G_V(x, y) \le  G_{\overline{B}(z_0, r)^c} (x, y)   \le c_1  \frac{\Phi(r)}{V(r)}= c_1  \frac{\Phi(r)}{V(r)} \phi(x),
\end{align*}
with $c_1=c_1(b)$.

For the remainder of the proof, we assume that  $x \in U:=V \cap B(z_0, ar)$. Let $f \ge0$ be supported on $A(z_0, 2r, 4r)$ with $\int f(w) m(dw)=1$.
Then, by the strong Markov property,
$$
G^\psi f(x) =\pi^\psi_U (G^\psi f)(x)=\pi^\psi_U ( {\bf 1}_{\overline{A}(z_0, ar, br)}G^\psi f)(x)+\pi^\psi_U ( {\bf 1}_{B(z_0, br)^c} G^\psi f)(x)=:I+II.
$$
First note that by  {\bf D2}$(z_0, r_0)$ and \eqref{e:Gmon}, for $y \in \overline{A}(z_0, ar, br)$,
$$
G^\psi f(y) \le \int_{\overline{A}(z_0, 2r, 4r)} G_V(y,w) f(w)m(dw)  \le \int_{\overline{A}(z_0, 2r, 4r)} G_{\overline{B}(z_0, r)^c}(y,w) f(w)m(dw) \le c_2 \frac{\Phi(r)}{V(r)}.
$$
Thus, combining  this with Lemma \ref{l:4.4} we get
\begin{align*}
&I \le \left(\sup_{y \in \overline{A}(z_0, ar, br)} G^\psi f(y) \right)
  \pi^\psi_U(x, \overline{A}(z_0, ar, br)) \le   c_2\frac{\Phi(r)}{V(r)}  \pi^\psi_U(x, V \setminus U)  \le c_2  \frac{\Phi(r)}{V(r)} \phi(x).
\end{align*}
For $II$, by using \eqref{e:AC2} and  Lemma \ref{l:etauU}, we get that for $z \in \overline{B}(z_0, br)^c$,
\begin{eqnarray*}
  \int_U G^\psi_U(x, y) j(y, z) m(dy) &\le& c_3 \int_U G^\psi_U(x, y) m(dy)  j(z_0, z) \\
  &\le &c_4
  \phi(x) \Phi(r)  j(z_0, z),
\end{eqnarray*}
  with $c_4=c_4(b/a)$.
  Thus  by \eqref{e:Gmon} and \eqref{e:psilevy},
\begin{eqnarray*}
II &\le &c_4
 \phi(x) \Phi(r)  \int_{ \overline{B}(z_0, br)^c }
G^\psi f(z)j(z_0, z) m(dz) \\
&\le &
c_4
 \phi(x) \Phi(r)
 \int_{A(z_0, 2r, 4r)}   \int_{ \overline{B}(z_0, br)^c } G_V(z,y)  j(z_0, z)m(dz)   f(y)m(dy)\\
 &\le &
c _4
 \phi(x) \Phi(r)
 \int_{A(z_0, 2r, 4r)}  \widehat{P}_{\overline{B}(z_0, r)^c }  (y, z_0)   f(y)m(dy).
\end{eqnarray*}
Finally using
the dual version of Proposition \ref{p:pcom}, we conclude that
\begin{align*}
 II \le c_5\phi(x) \int_{A(z_0, 2r, 4r)}     f(y)m(dy)   \frac{\Phi(r)}{V(r)} = c_5  \frac{\Phi(r)}{V(r)} \phi(x)\, .
\end{align*}
\qed

\begin{lemma}\label{l:4.9}
There exists a constant $c=c(a,z_0)>0$ such that for all $x \in \overline{A}(z_0, 2r, 4r)$,
\begin{align}
\label{e:pipsiest}
\pi^\psi(x, dy)/m(dy) \le  \frac{c}{V(r)}  \ind_{\overline{B}(z_0, ar)}(y).
\end{align}
\end{lemma}
\pf
Let $b:=a/2+1$ so that $b\in (a,2)$. First note that $\psi$ vanishes on $\X \setminus B(z_0, ar)$.
Thus
\begin{align} \label{e:Vsupp}
\pi^\psi(y, \X \setminus B(z_0, ar)) = 0, \quad y \in V.
\end{align}
Fix $x \in \overline{A}(z_0, 2r, 4r)$. If $f$ is a non-negative function on $\X$ vanishing in $\X \setminus (B(z_0, ar) \cap V)$, then by \eqref{e:pipsi1} and the dual version of Lemma \ref{l:4.5} (together with $G^{\psi}(x,y)=\widehat{G}^{\psi}(y,x)$),
\begin{align}\label{e:pipsi11}
 \pi^\psi f(x)  \le   c_1  \frac{\Phi(r)}{V(r)} \int_{V \cap B(z_0, ar)} \phi(y)    \psi(y) f(y) m(dy).
\end{align}
Since for $y\in B(z_0,ar)$ we have $\phi(y)    \psi(y)  \le c_2(\Phi(r))^{-1}$ by the definition in
\eqref{e:defpsi} and assumption {\bf B2-a}$(z_0, r_0)$ (with a constant $c_3=c_3(a)$ possibly different from the one in
{\bf B2-a}$(z_0, r_0)$), we have
\begin{align}\label{e:pipsi12}
\pi^\psi f(x)  \le    \frac{c_3}{V(r)}\int_{V \cap B(z_0, ar)} f(y) m(dy).
\end{align}
On the other hand, if $g \in \dom(\sA)$ vanishes in $V$ then by \eqref{e:pipsi2},
\begin{align}\label{e:pipsi21}
\pi^\psi g(x)   = \int_{\X \setminus V} \expr{\int_V G^\psi(x, y) j(y, z) m(dy)} g(z) m(dz).
\end{align}
Assume $z \in \X \setminus V \subset B(z_0, ar)$ and let
$$
I:= \int_{V \cap  B(z_0, br)}  G^\psi(x, y) j(y, z) m(dy)\quad \text{and}
\quad II:=\int_{\overline{B}(z_0, br)^c} G^\psi(x, y) j(y, z) m(dy).
$$
We now consider $I$ and $II$ separately.

By the dual version of Lemma \ref{l:4.5}, assumption {\bf B2-a}$(z_0, r_0)$ and
the fact that $\widehat{\sA}\phi(z)=\widehat{J}\phi(z)$, for some  $c_4=c_4(a)>0$,
\begin{align}\label{e:pipsi23}
I \le c_4  \frac{\Phi(r)}{V(r)}   \int_{V \cap  B(z_0, br)}  \phi(y)    j(y, z) m(dy)
= c_4  \frac{\Phi(r)}{V(r)}  \widehat  \sA \phi(z) \le  \frac{c_4}{V(r)}.
\end{align}

On the other hand, by assumption {\bf C2}$(z_0, r_0)$ and \eqref{e:Gmon}, for some $c_5=c_5(a)>0$,
\begin{align}\label{e:pipsi24}
II &\le c_5 \int_{\overline{B}(z_0, br)^c} G^\psi(x, y) j(y, z_0) m(dy)\nonumber\\
 &\le c_5 \int_{\overline{B}(z_0, br)^c} G_V(x, y) j(y, z_0) m(dy)\le c_5P_V(x, z_0),
\end{align}
which is less than or equal to  $c_6 V(r)^{-1}$ by Proposition \ref{p:pcom}.
Therefore, \eqref{e:pipsi23} and \eqref{e:pipsi24} imply that for all $g \in \dom(\sA)$ vanishing in $V$ we have
\begin{align}\label{e:pipsi25}
 \pi^\psi g(x)   \le  \frac{c_7}{V(r)} \int_{\X \setminus V} g(z) m(dz).
\end{align}
Since $\dom(\sA)$ is dense in $C_0(\sX)$ we have $\pi_\psi(x, dy)/m(dy) \le  {c_7}/{V(r)}$ on $V^c$ too.
\qed

\begin{corollary}\label{c:4.9}
Let $f$ be a non-negative function on $\X$ and $x$ a point in $\overline{A}(z_0, 2r, 4r)$ such that $f(x)\le \E_x f(X_{\tau})$ for every stopping time $\tau\le \tau_{\overline{B}(z_0, r)^c}$. Then
\begin{equation}\label{e:har:c}
f(x) \le \frac{c}{V(r)} \int_{\overline{B}(z_0, ar)} f(y) m(dy)\, ,
\end{equation}
where $c=c(a)$ is the constant from Lemma \ref{l:4.9}.
\end{corollary}
\pf
Recall from \eqref{e:psi-tau-a} that $\pi^\psi f(x)=\int_0^\infty \E_x(f(X_{\tau_a})) e^{-a} da$. Since $\tau_a\le \tau_V\le \tau_{\overline{B}(z_0,r)^c}$, we have that $f(x) \le \E_x f(X_{\tau_a})$, and therefore $f(x) \le \pi^\psi f(x)$. Thus by \eqref{e:pipsiest},
$$
f(x) \le \int f(y) \pi^\psi(x, dy) \le  \frac{c}{V(r)} \int_{\overline{B}(z_0, ar)} f(y) m(dy).
$$
\qed

\begin{lemma}\label{l:exit}
There exists $c=c(z_0, a)>0$ such that for any $r \ge r_0$ and any open set $D\subset \overline{B}(z_0, r)^c$ we have
$$
 {\P}_x\left(X_{\tau_D} \in B(z_0, r)\right) \,\le\, c\,V(r) P_D(x,z_0)\, , \qquad x \in D\cap \overline{B}(z_0, ar)^c\, .
$$
\end{lemma}
\pf
\eqref{e:Assumptions B2(b)} says that for any $f\in \dom(\overline{B}(z_0, r), B(z_0, ar))$, we have
\begin{align*}
 \max(\sA f(z), \hat{\sA} f(z)) \le   c {\bf 1}_{B(z_0, r)^c} (z)  V(r) j(z,z_0)
\end{align*}
for some $c(z_0, a)>0$ independent of $r \ge r_0$.
Thus, by Dynkin's formula we have
\begin{align}
& \E_x\left[ f(X_{\tau_D})\right] = \int_D G_D(x,z)  {\sA}  f(z)\, dz\nonumber \\
& \le  c V(r) \int_D G_D(x,z)  j(z,z_0)  \, dz= c V(r)P_D(x,z_0).
\label{e:dynkin}
\end{align}
Finally, since ${\bf 1}_{\overline{B}(z_0,r)}\le f$,
$\P_x(X_{\tau_D}\in \overline{B}(z_0,r))\le \E_x\left[ f(X_{\tau_D})\right] \le c V(z_0, r) P_D(x,z_0)$.
\qed

\begin{prop}\label{p:main-aux}
Let $b\in (a,2)$. There exists $c=c(z_0,a,b)>1$ such that for any $r \ge r_0$, any open
set $D \subset \overline{B}(z_0,r)^c$ and any non-negative function $u$ on $\X$ which
is regular harmonic with respect to $X$ in $D$ and
vanishes on $\overline{B}(z_0,r)^c\cap\left(\overline{D}^c\cup D^{\mathrm{reg}}\right)$, it holds that
\begin{align}\label{e:bhp-inf-lemma}
&c^{-1}P_{D \cap \overline{B}(z_0, 2br)^c}(x,z_0) \int_{B(z_0, 2ar)}u(z)\, m(dz) \le u(x) \nonumber\\
 &\,\,\,\, \le c P_{D \cap \overline{B}(z_0, 2br)^c}(x,z_0) \int_{B(z_0,2ar)}u(z)\, m(dz)
\end{align}
for all $x\in D\cap \overline{B}(z_0,4r)^c$.
\end{prop}
\pf
Let $O:=D \cap \overline{B}(z_0, 2br)^c$, $D_1:=\overline{A}(z_0, 2ar, 2br)$ and $D_2:=B(z_0, 2ar)$.
By the harmonicity of $u$,
\begin{align}\label{e:bhp1}
 u(x)=\E_x[u(X_{\tau_{O}})]=\E_x[u(X_{\tau_{O}}):X_{\tau_{O}}\in D_1 ]+\E_x[u(X_{\tau_{O}}):X_{\tau_{O}}\in D_2 ]\, , \ \ x\in D.
\end{align}
Let $D^{\mathrm{irr}}$ be the set of points in $\partial D$ which are irregular for $D^c$ with respect to $X$.
Since $u$ vanishes on $\overline{B}(z_0,r)^c\cap\left(\overline{D}^c\cup D^{\mathrm{reg}}\right)$,
it follows that $u(y)\le \E_y u(X_{\tau})$ for every stopping time $\tau\le \tau_{\overline{B}(z_0,r)^c}$
and every $y\in \overline{B}(z_0,r)^c\setminus D^{\mathrm{irr}}$.
Since $D^{\mathrm{irr}}$ is polar with respect to $X$, we see that $X_{\tau_O}\notin D^{\mathrm{irr}}$.
It  follows from Corollary \ref{c:4.9}
and Lemma \ref{l:exit} that for all $x\in D\cap \overline{B}(z_0,4r)^c$,
\begin{align}\label{e:bhp2}
&\E_x[u(X_{\tau_{O}}):X_{\tau_{O}}\in D_1 ]
\le \left(\sup_{y \in D_1\setminus D^{\mathrm{irr}}} u(y) \right) \P_x(X_{\tau_{O}}\in D_1) \nonumber\\
&\le c_1 V(br)P_O(x,z_0) \frac{c_2}{V(r)} \int_{B(z_0,a^2r)}u(z)\, m(dz) \\
&\le c_3 P_{O} (x,z_0)  \int_{B(z_0,2ar)}u(z)\, m(dz),
\end{align}
where $c_3=c_3(a,b)$.
On the other hand, by assumption {\bf C2}$(z_0, r_0)$, for all $x\in D\cap \overline{B}(z_0,4r)^c$,
\begin{align}\label{e:bhp3}
\E_x[u(X_{\tau_{O}}):X_{\tau_{O}}\in D_2]
&=\int_{B(z_0, 2ar)} \int_{O} G_{O}(x,y)j(y,z)m(dy)u(z)m(dz)\nonumber\\
&\asymp
\int_{B(z_0, 2ar)} \int_{O} G_{O}(x,y)j(y,z_0)m(dy)u(z)m(dz)\nonumber\\
&=P_{O}(x,z_0)\int_{B(z_0, 2ar)}u(z)m(dz).
\end{align}
The proposition now follows from \eqref{e:bhp1}--\eqref{e:bhp3}.
\qed

\begin{lemma}\label{l:poisson-comparison}
For any $b\in (a,2)$ there exists $c=c(a,b)>0$ such that for every
$r\ge r_0$ and every open set $D\subset \overline{B}(z_0,r)^c$,
$$
P_{D\cap \overline{B}(z_0,2br)^c}(x,z_0)\le P_D(x,z_0)\le c P_{D\cap \overline{B}(z_0,2br)^c}(x,z_0)\,  , \qquad x\in D\cap \overline{B}(z_0,2ab r)^c\, .
$$
\end{lemma}
\pf First note that, since $D$ is Greenian, by the strong Markov property
for all open set $U \subset D$,
$G_D(x,y)= G_U(x,y) + \E_x\left[   G_D(X_{\tau_U}, y); \tau_U<\infty\right]$ for every $(x,y) \in \X \times \X$.
Thus
\begin{align*}
&P_{D}   (x,z_0) =P_{D\cap \overline{B}(z_0, 2br)^c}   (x,z_0)\\
&\quad + \E_x [P_{D}   (X_{\tau_{D\cap \overline{B}(z_0, 2br)^c}}, z_0) : X_{\tau_{D\cap \overline{B}(z_0, 2b r)^c}} \in  {B}(z_0, 2br) \setminus \overline{B}(z_0, r), \tau_{D\cap \overline{B}(z_0, 2b r)^c}<\infty].
\end{align*}
By Proposition \ref{p:pcom}, Lemma \ref{l:exit} and the doubling property,
for $x\in D\cap \overline{B}(z_0, 2ab r)^c$,
\begin{eqnarray*}
\lefteqn{ \E_x [P_D  (X_{\tau_{D\cap \overline{B}(z_0, 2br)^c}}, z_0):\,  X_{\tau_{D\cap \overline{B}(z_0, 2br)^c}} \in  B(z_0, 2br) \setminus \overline{B}(z_0, r)]}\\
&\le & \left(\sup_{z \in {B}(z_0, 2br) \setminus \overline{B}(z_0, r)} P_{D}   (z ,z_0) \right)
\P_x (X_{\tau_{D\cap \overline{B}(z_0, 2br)^c}} \in  \overline{B}(z_0, 2br))\\
&\le & c_1 \frac{V(2br)}{V(r)} P_{D\cap \overline{B}(z_0,2br)^c}
\le  c_2 P_{D\cap \overline{B}(z_0, 2br)^c}   (x,z_0)\, .
\end{eqnarray*}
This finishes the proof. \qed

\noindent
{\bf Proof of Theorem \ref{t:main-infty}} Let $a\in (1,2)$ and choose $b=a/2+1$.
Let $D\subset \overline{B}(z_0,r)^c$ and let $u$ be a non-negative function on $\X$
which is regular harmonic with respect to $X$ in $D$ and
vanishes on $\overline{B}(z_0,r)^c\cap\left(\overline{D}^c\cup D^{\mathrm{reg}}\right)$.
Since $\overline{B}(z_0, 8r)^c\subset \overline{B}(z_0, 4r)^c\cap \overline{B}(z_0, 2abr)^c$,
it follows from Proposition \ref{p:main-aux} and Lemma \ref{l:poisson-comparison} that
$$
u(x)\asymp P_D(x,z_0)\int_{B(z_0,2ar)} u(y)m(dy)\, ,\qquad x\in  D\cap \overline{B}(z_0, 8r)^c\,  ,
$$
with a constant depending on $a$.
\qed

\section{Finite boundary point  }\label{s:finite}
The goal of this section is to state an analog of Theorem \ref{t:main-infty} for finite boundary points. Again, recall that assumptions {\bf A},  {\bf B} and  {\bf C} are in force.
Recall that $R_0\in (0,\infty]$ is the localization radius and that $\sA$ and $\widehat \sA$ are the generators of
$(P_t)$ and $(\widehat{P}_t)$ in $C_0(\X)$.
The processes $X$ and $\widehat{X}$ are not assumed to be conservative.

Similarly as in Section \ref{s:main-result} we fix a point $z_0\in \X$
which now serves as a boundary point of an open set.
For $r>0$, we let $V(r)=V(z_0,r):=m(B(z_0,r))$ and assume that the volume function $V:[0,\infty)\to [0,\infty)$ satisfies \eqref{e:V-doubling}
and, instead of \eqref{e:V-another},
we assume that
there exist $c>1$ and $r_0 \in (0, R_0]$ and $n_0 \in \N$ with $n_0 \ge 2$ such that
\begin{equation}\label{e:V-another1}
V(n_0r) \ge c V(r), \quad r \le r_0\,  .
\end{equation}

We also assume the existence of an increasing function $\Phi=\Phi(z_0, \cdot):[0,\infty)\to [0,\infty)$ satisfying the doubling property \eqref{e:Phi-doubling} which again will be crucial in obtaining the scale invariant results.

Similarly as in Section \ref{s:main-result} we introduce some additional assumptions.
Recall the notation $\dom(K,D)$ from assumption {\bf B}. In the following assumptions, $r_0$ is
a number in $(0, R_0]$.

\smallskip
\noindent
{\bf B1-a}$(z_0, r_0)$: For any $a\in (1/2,1)$, there exists $c=c(z_0,a)$ such that for any $r<r_0$,
$$
 \rho(r) := \inf_{f\in \dom(\overline{B}(z_0, ar), B(z_0, r)) } \sup_{x\in \X} \max(\sA f(x), \widehat{\sA} f(x))\le  \frac{c}{\Phi(r)}.
$$

\noindent
{\bf B1-b}$(z_0, r_0)$: For any $a\in (1/2,1)$, there exists $c=c(z_0,,a)$ such that for any $r<r_0$
and any $f\in \dom(\overline{B}(z_0, ar), B(z_0, r))$,
\begin{align*}
 \max(\sA f(x), \widehat{\sA} f(x)) \le  c V(r) j(x,z_0),  \quad x \in \overline{A}(z_0, ar, (a+1)r).
\end{align*}

\noindent
{\bf B1-c}$(z_0, r_0)$: For any $1/2<b<a<1$, there exists $c=c(z_0,,a,b)$ such that for any $r<r_0$,
$$
\inf_{f\in \dom(\overline{A}(z_0, br, ar), A(z_0, r/2, r)) } \sup_{x \in \X} \max(\sA f(x), \widehat{\sA} f(x))\le  \frac{c}{\Phi(r)}.
$$

To assumption {\bf C} we add

\smallskip
\noindent
{\bf C1}$(z_0, r_0)$: For any $a\in (1/2,1)$,  there exists $c= c(a, r_0, z_0)$ such that for $r <r_0$,
$x \in B(z_0, ar)$ and $y \in \X \setminus B(z_0, r)$,
\begin{align}\label{e:AC2-finite}
 c^{-1}{j(z_0, y)} \le {j(x, y)} \le c{j(z_0, y)} , &&
 c^{-1}\, {\widehat j(z_0, y)} \le {\widehat j(x, y)} \le c{\widehat j(z_0, y)},
\end{align}
and
\begin{equation}\label{e:AC21-finite}
\inf_{y\in \overline{A}(z_0,ar,r)} \min(j(z_0,y), \widehat{j}(z_0,y))\ge \frac{c}{V(r)\Phi(r)}\, .
\end{equation}

Note that by assumption {\bf C1}$(z_0, r_0)$,  the function $f$ in assumption {\bf B1-b}$(z_0, r_0)$ satisfies
\begin{align}\label{e:Assumptions B1(b)}
 \max(\sA f(x), \widehat{\sA} f(x)) \le c  {\bf 1}_{B(z_0, ar)^c} (x)  V(r) j(x,z_0)
\end{align}
for a constant $c=c(z_0,a)>0$.
In fact,  for $f \in  \dom(\sA) \cap \dom(\widehat{\sA})$ such that $f(x) = 1$ for $x \in   B(z_0, ar)$, $f(x) = 0$ for $x \in \X \setminus \overline{B}(z_0, r)$ and $0 \le f(x) \le 1$ for $x \in \X$, we have
\begin{align*}
\sA f(x) {\bf 1}_{\X \setminus \overline{B}(z_0, (a+1)r)} (x) & =
{\bf 1}_{\X \setminus \overline{B}(z_0, (a+1)r)} (x)\int_{B(z_0, r)} f(y) j(x,y)m(dy) \nonumber\\
&\le c {\bf 1}_{\X \setminus \overline{B}(z_0, (a+1)r)} (x) V(r)j(x,z_0).
\end{align*}

The final assumption concerns Green functions  of balls.

\smallskip
\noindent
{\bf D1}$(z_0, r_0)$: For any $a\in (1/2,1)$ there exists a constant $c=c(z_0, r_0, a)$ such that for all $r<r_0$,
 $$
\sup_{x \in \overline{A}(z_0, ar,r)} \sup_{y \in  B(z_0,r/2)} \max(G_{B(z_0, r)}(x, y), \widehat{G}_{B(z_0, r)}(x, y)) \le c\frac{\Phi(r)}{V(r)}.
$$
\smallskip

\begin{thm}\label{t:main-finite}
Let $z_0\in \X$.
Assume that, in addition to {\bf A}, {\bf B} and {\bf C}, the assumptions \eqref{e:V-doubling}, \eqref{e:V-another1},  \eqref{e:Phi-doubling}, {\bf B1-a}$(z_0, r_0)$,
{\bf B1-b}$(z_0, r_0)$, {\bf B1-c}$(z_0, r_0)$, {\bf C1}$(z_0, r_0)$ and {\bf D1}$(z_0, r_0)$
hold true for some $r_0\in (0, R_0]$.
For any $a\in (1/2,1)$, there exists $C_1=C_1(z_0, r_0, a)>1$ such that for any
$r < r_0/(2n_0)$,
any open set $D \subset B(z_0,r)$ and any non-negative function $u$ on $\X$ which is
regular harmonic with respect to $X$ in $D$ and
vanishes on $B(z_0,r)\cap \left(\overline{D}^c\cup D^{\mathrm{reg}} \right)$, it holds that
\begin{align}\label{e:bhp-fin}
 &C_1^{-1} \, \E_x \tau_{D}\int_{\overline{B}(z_0,ar/2)^c} j(z_0,z)u(z)\, m(dz) \le u(x)\nonumber\\
&\,\,\,\, \le C_1 \, \E_x\tau_{D}\int_{\overline{B}(z_0, ar/2)^c} j(z_0,z)u(z)\, m(dz)
\end{align}
for all $x\in D\cap B(z_0,r/8)$.
\end{thm}

As a consequence of Theorem \ref{t:main-infty}, one immediately gets the following scale invariant
uniform boundary Harnack principle.

\begin{corollary}[Boundary Harnack Principle]\label{c:bhp-finite}
Let $z_0\in \X$.
Suppose that, in addition to {\bf A}, {\bf B} and {\bf C}, the assumptions \eqref{e:V-doubling}, \eqref{e:V-another1},  \eqref{e:Phi-doubling}, {\bf B1-a}$(z_0, r_0)$,
{\bf B1-b}$(z_0, r_0)$, {\bf B1-c}$(z_0, r_0)$, {\bf C1}$(z_0, r_0)$ and {\bf D1}$(z_0, r_0)$
hold true for some $r_0\in (0, R_0]$.
There exists $C_2=C_2(z_0, r_0)>1$ such that for any
$r < r_0/(2n_0)$, any open set $D \subset B(z_0,r)$
and any non-negative
functions $u$ and $v$ on $\X$ which are regular harmonic with respect to $X$ in $D$ and vanish
on $B(z_0,r)\cap \left(\overline{D}^c\cup D^{\mathrm{reg}} \right)$, it holds that
\begin{equation}\label{e:bhp-fin-cor}
    C_2^{-1} \frac{u(y)}{v(y)} \le  \frac{u(x)}{v(x)}\le C_2 \frac{u(y)}{v(y)}\, ,\qquad \textrm{for all }x,y\in D\cap B(z_0,r/8)\,  .
\end{equation}
\end{corollary}

\section{Examples}\label{s:examples}

\begin{example}\label{ex:umlp}
{\rm
Let $X=(X_t, \P_x)$ be a purely discontinuous
symmetric L\'evy process in $\R^d$ with L\'evy exponent $\Psi(\xi)$ so that
$$
\E_x\left[e^{i\xi\cdot(X_t-z_0)}\right]=e^{-t\Psi(\xi)}, \qquad t>0, x\in \R^d, \xi\in\R^d.
$$
Thus the state space $\X=\R^d$, the measure $m$ is the $d$-dimensional Lebesgue measure and the
localization radius $R_0=\infty$.
Assume that $r\mapsto j_0(r)$ is a strictly positive and
nonincreasing function on $(0, \infty)$ satisfying
\begin{equation}\label{e:fuku1.1}
j_0(r)\le cj_0(r+1), \qquad r>1,
\end{equation}
for some $c>1$, and the L\'evy measure of $X$ has a density $J$ such that
\begin{equation}\label{e:fuku1.2}
\gamma^{-1}j_0(|y|)\le J(y) \le \gamma j_0(|y|), \qquad y\in \R^d,
\end{equation}
for some $\gamma>1$.
Since
$\int_0^\infty j_0(r) (1\wedge r^2) r^{d-1}dr < \infty$ by \eqref{e:fuku1.2},
the function $x \to j_0(|x|)$ is the L\'evy density of  an isotropic unimodal L\'evy process whose characteristic exponent is
\begin{equation}\label{e:fuku1.3}
\Psi_0(|\xi|)= \int_{\R^d}(1-\cos(\xi\cdot y))j_0(|y|)dy.
\end{equation}
The L\'evy exponent $\Psi$ can be written as
$$
\Psi(\xi)= \int_{\R^d}(1-\cos(\xi\cdot y))J(y)dy
$$
and, clearly by \eqref{e:fuku1.2}, it satisfies
\begin{equation}\label{e:fuku1.4}
\gamma^{-1} \Psi_0(|\xi|)\le \Psi(\xi) \le \gamma \Psi_0(|\xi|),
\quad \mbox{for all } \xi\in \R^d\, .
\end{equation}
 The function $\Psi_0$ may not be increasing.
However, if we put $\Psi^*_0(r):= \sup_{s \le r} \Psi_0(s)$, then,
by \cite[Proposition 2]{BGR14}  (cf.~also \cite[Proposition 1]{G}), we have
$$
\Psi_0(r) \le\Psi^*_0(r) \le \pi^2 \Psi_0(r).
$$
Thus by \eqref{e:fuku1.4},
\begin{equation}\label{e:fuku1.5}
(\pi^2\gamma)^{-1}  \Psi^*_0(|\xi|)\le \Psi(\xi) \le \gamma \Psi^*_0(|\xi|),
\quad \mbox{for all } \xi  \in \R^d\, .
\end{equation}
Under the above assumptions, our process $X$ obviously satisfies Assumptions
{\bf A}, {\bf B}  and {\bf C}.

Let $\Phi(r)=(\Psi^*_0(r^{-1}))^{-1}$.
Since $X$ is a purely discontinuous symmetric L\'evy process, we can write down the generator $\sA$ of $X$
explicitly in terms of the L\'evy density.
Using this explicit formula and \cite[Corollary 1]{G}
one can easily check that
Assumptions {\bf B1-a}$(z_0, r_0)$, {\bf B1-b}$(z_0, r_0)$, {\bf B1-c}$(z_0, r_0)$ and {\bf B2-a}$(z_0, r_0)$ are also satisfied
for all $z_0\in \R^d$
and $r_0>0$.

Suppose now that $\Psi_0$ satisfies the following scaling condition at infinity:

\noindent
{\bf H1}:
There exist constants $0<\delta_1\le \delta_2 <1$ and $a_1, a_2>0$  such that
\begin{equation}\label{e:fuku1.6}
a_1\left(\frac{t}{s}\right)^{2\delta_1} \le \frac{\Psi_0( t)}{\Psi_0( s)} \le a_2 \left(\frac{t}{s}\right)^{2\delta_2} , \quad t\ge s \ge 1\, .
\end{equation}
Then by  \cite[(15) and Corollary 22]{BGR14},
for every $R>0$, there exists $c=c(R)>1$ such that
\begin{equation}\label{e:fuku1.7}
c^{-1}\frac{\Psi_0(r^{-1})}{r^d} \le j_0(r)\le c \frac{\Psi_0(r^{-1})}{r^d}
\quad \hbox{for } r\in (0, R].
\end{equation}
Using \eqref{e:fuku1.1}, \eqref{e:fuku1.7}
and \cite[Lemma 2.7]{KSV15}, one can easily see that, there exists $r_1>0$ such that
Assumption {\bf C1}$(z_0, r_0)$ and Assumption {\bf D1}$(z_0, r_0)$ are satisfied for all $z_0\in \R^d$
and $r_0\le (0, r_1]$.

Now we assume, instead of  {\bf H1}, that $\Psi_0$ satisfies the following scaling condition at the origin:

\noindent
{\bf H2}:
There exist constants $0<\delta_3\le \delta_4 <1$ and $a_3, a_4>0$  such that
\begin{equation}\label{e:fuku1.6at0}
a_3\left(\frac{t}{s}\right)^{2\delta_3} \le \frac{\Psi_0( t)}{\Psi_0( s)} \le a_4 \left(\frac{t}{s}\right)^{2\delta_4} , \quad s\le t \le 1\, .
\end{equation}
It follows from
\cite[Corollary 7]{BGR14} or \cite[Theorem 2.2]{KKK}
 (see \cite[(15)]{BGR14}) and \eqref{e:fuku1.2} that, there exists $c_1>1$ such that
\begin{equation}\label{e:fuku1.7gl}
J(x) \le \gamma j_0(|x|)\le c \gamma \frac{\Psi_0(|x|^{-1})}{r^d}
\quad \hbox{for } x \in \R^d\setminus\{0\}.
\end{equation}

We now prove a matching lower bound for $j_0$ away from the origin. The proof is similar to that of \cite{BGR14}.

Let $Y$ be the isotropic unimodal L\'evy process whose characteristic exponent is
$\Psi_0(|\xi|)$ and  $x \to p^0_t(|x|)$ be its transition density.
Let
$
{f}_t(r):=\P(|Y_t|^2>r)$ for $r \ge 0$ and $ t>0$.
Then, using  \cite[Lemma 4 and (13)]{BGR14}, ${\cal L}{{f}_t}$,  the Laplace transform of $f_t$,
satisfies that for all $0<u<v \le 1$,
$$\frac{\mathcal{L}{f}_t(v )}{\mathcal{L}{f}_t(u)}\le c_2 (v/u)^{-1}\frac{1-e^{-\pi^2t\Psi_0(\sqrt{v})}}{1-e^{-t\Psi_0(\sqrt{u})}}\leq c_2 (v/u)^{-1}\frac{1-e^{-\pi^2t\Psi_0( \sqrt{u})a_4 (v/u)^{\delta_4}}}{1-e^{-t\Psi_0(\sqrt{u})}}\leq c_3(v/u)^{\delta_4-1}.$$
Thus, by \cite[Proposition 2.3]{KSV14a} and \cite[Lemma 4]{BGR14},
\begin{equation}\label{e:tail0}
\P(|Y_t|\ge r)=f_t(r^{2})\geq c_4 \mathcal{L}f_t(
r^{-2})\geq 2c_5
(1-e^{-t\Psi_0^*(
1/r
)}),\qquad r \ge 1.
\end{equation}
Let
$a\ge 2$. Since $r \to p^0_t(r)$ is decreasing,
 we have
\begin{equation}\label{e:tail1}p^0_t(r)\ge
\frac{
{\P}
 (  r\leq |Y_t|<a r)}{ \left|  B(0, a r)\setminus B(0, r)\right|}
= \frac {c_6} {a^d-1}  r^{-d}(P (   |Y_t|\ge r) - \P (   |Y_t|\ge a r)).\end{equation}
Let $r \ge 1$ and  $t\Psi_0^*(1/r)\leq 1$.
Using \eqref{e:tail0}, the inequality $s/2\le 1-e^{-s}\le s$ for $s\in (0, 1]$,
and  \cite[Corollary 6]{BGR14}, we get
\begin{align}
\label{e:tail2}
\P (   |Y_t|\ge r) - \P (   |Y_t|\ge a r)
&\ge c_5t\Psi_0^*(1/r)-
\frac{2e}{e-1}(2d+1)
t\Psi_0^*(1/ar) \nn\\
&\ge c_5t\Psi_0^*(1/r)
(1-c_7\frac{\Psi_0^*(1/ar)} {\Psi_0^*(1/r)}).
\end{align}
Choose $a\ge 2$ large enough so that
for $ar \ge 1$,
\begin{align}
\label{e:tail3}
c_7\frac{\Psi_0^*(1/ar)} {\Psi_0^*(1/r)}\le
c_7 a_3^{-1} \kappa^{-2\delta_3} \le \frac12.
\end{align}
Then, combining \eqref{e:tail1}--\eqref{e:tail3}, we obtain
$$
p^0_t(r)\geq c_8 t \Psi_0^* (1/r) r^{-d},    \quad r \ge 1/a, t\Psi_0^*(1/r)\leq 1,
$$
which, together with \eqref{e:fuku1.2} and the fact that $J$ is the weak limit of   $p^0_t$,  implies
\begin{equation}\label{e:fuku1.7gln}
 J(x)  \ge \gamma^{-1} j_0(|x|)\ge c_9 \frac{\Psi_0(|x|^{-1})}{|x|^d}
\quad \hbox{for } |x|\ge 1.
\end{equation}
Hence, {\bf C2}$(z_0,r_0)$ is valid.

If $d \ge 3$, by \cite[(14) and p. 26]{G}, the Green function $G(x)$ of $X$
has the following upper bound:
$$
\int_{B(0, r)}G(x) dx \le  \frac{c_{10}}{\Psi_0(r^{-1})}.
$$
We further assume that there exists a positive constant $c_{11}>1$ and a non-increasing function $r \to G_0(r)$ such that
\begin{equation}\label{e:2sgfngle}
 c^{-1}_{11} G_0(|x|) \le G(x) \le c_{11} G_0(|x|), \quad x \in \R^d.
\end{equation}
Then we have that for all $x \in \R^d$,
\begin{eqnarray}
G(x)&\le&    c_{11} G_0(|x|)  \le
c_{12}|x|^{-d} \int_{B(0, |x|)}G_0(|y|) dy\nonumber\\
& \le& c_{12} c_{11} |x|^{-d} \int_{B(0, |x|)}G(y) dy  \le
 \frac{c_{12}c_{13}c_{10}}{|x|^d\Psi_0(|x|^{-1})}.\label{e:upgfngle}
\end{eqnarray}

It follows immediately from \eqref{e:fuku1.7gl}, \eqref{e:fuku1.7gln} and \eqref{e:upgfngle} (when $d \ge 3$) that {\bf C2}$(z_0, r_0)$
and {\bf D2}$(z_0, r_0)$ are satisfied
for all $z_0\in\R^d$ and all $r_0\ge 1$.
We remark here that, by \cite[Lemma 3.3]{KSV12}, the upper bound
$
G(x) \le
 \frac{c}{|x|^d\Psi_0(|x|^{-1})}$
 holds for $d > 2\delta_4$
when $X$
 is a subordinate Brownian motion
whose
Laplace exponent $\phi$ is a complete Bernstein function
and that  $\xi \mapsto\phi(|\xi|^2)$ satisfies Assumption {\bf H2}.

Using Assumption {\bf H2} and the explicit form of the
the generator, one can easily check that Assumption {\bf B2-b}$(z_0, r_0)$
is also satisfied for all $z_0\in\R^d$ and all $r_0 \ge 1$ (e.g.~see \cite[(3.4)]{KSV14}).
Thus, under the assumptions above, {\bf B2-a}$(z_0, r_0)$, {\bf B2-b}$(z_0, r_0)$,
{\bf C2}$(z_0, r_0)$ and {\bf D2}$(z_0, r_0)$ all hold.

Suppose that $\alpha\in (0, 2)$. The subordinate Brownian motion in $\R^d$ via a subordinator
with Laplace exponent $\phi(\lambda)=\log(1+\lambda^{\alpha/2})$ is called a geometric
$\alpha$-stable process. Let $\phi^{(1)}(\lambda)=\log(1+\lambda^{\alpha/2})$. For $n>1$, let
$\phi^{(n)}(\lambda)=\phi^{(1)}(\phi^{(n-1)}(\lambda))$. A subordinate Brownian motion in $\R^d$ via
a subordinator with Laplace exponent $\phi^{(n)}$ is called an $n$-iterated geometric $\alpha$-stable
process. It is clear that geometric $\alpha$-stable and $n$-iterated geometric $\alpha$-stable
processes satisfy condition {\bf H2} and \eqref{e:2sgfngle} and, again by \cite[Lemma 3.3]{KSV12},    for $d > 2\alpha$  the upper bound
$ G(x)\le  \frac{c_{(n)}}{|x|^d\phi^{(n)}(|x|^{-2})}$ holds.
Hence, for the geometric $\alpha$-stable and $n$-iterated geometric $\alpha$-stable
processes, {\bf B2-a}$(z_0, r_0)$, {\bf B2-b}$(z_0, r_0)$,
{\bf C2}$(z_0, r_0)$ and {\bf D2}$(z_0, r_0)$ all hold.
}
\end{example}

\begin{example}\label{ex:salphastable}
{\rm
Suppose that $\alpha\in (0, 2)$, $d\ge 2$ and that $X=(X_t, \P_x)$ is a strictly
$\alpha$-stable process in $\R^d$. Let $S$ be the unit sphere $S=\{x\in \R^d:
|x|=1\}$. For $\alpha\in (0, 1)$, $X$ is strictly $\alpha$-stable if and only if
$X$ is a L\'evy process with L\'evy exponent
$$
\Psi(\xi)=\int_S\lambda(d\theta)\int^\infty_0(1-e^{ir\theta\cdot\xi})r^{-1-\alpha}dr
$$
for some finite measure $\lambda$ on $S$. For $\alpha=1$,
$X$ is strictly $\alpha$-stable if and only if
$X$ is a L\'evy process with L\'evy exponent
$$
\Psi(\xi)=\int_S\lambda(d\theta)\int^\infty_0(1-e^{ir\theta\cdot\xi}+ir\theta\cdot\xi1_{(0, 1]})r^{-2}dr
+i\gamma\cdot \xi
$$
for some $\gamma\in\R^d$ and some finite measure $\lambda$ on $S$ satisfying $\int_S\theta\lambda(d\theta)=0$.
For $\alpha\in (1, 2)$,  $X$ is strictly $\alpha$-stable if and only if
$X$ is a L\'evy process with L\'evy exponent
$$
\Psi(\xi)=\int_S\lambda(d\theta)\int^\infty_0(1-e^{ir\theta\cdot\xi}+ir\theta\cdot\xi1_{(0, 1]})r^{-1-\alpha}dr
$$
for some finite measure $\lambda$ on $S$.

It follows from \cite{rao} that every semipolar set of $X$ is a polar set.
We will assume that $\lambda$ has a density with respect to the surface measure $\sigma$ on $S$ which
is bounded between two positive numbers.
Since $X$ is a a strictly $\alpha$-stable process, it automatically satisfies Assumption {\bf A}.
Let $\Phi(r)=r^{\alpha}$. By our assumption, it is obvious that Assumptions {\bf C},
and {\bf C1}$(z_0, r_0)$ and {\bf C2}$(z_0, r_0)$  are satisfied for all $z_0\in \R^d$ and $r_0\in (0, \infty)$.
It follows from \cite[(4.3)]{V} that
{\bf D1}$(z_0, r_0)$ and {\bf D2}$(z_0, r_0)$  are satisfied for all $z_0\in \R^d$ and $r_0\in (0, \infty)$. Since the generators of $X$ and its dual can be written out explicitly, one
can easily check that Assumptions {\bf B}, {\bf B1-a}$(z_0, r_0)$,  {\bf B1-b}$(z_0, r_0)$,
{\bf B1-c}$(z_0, r_0)$, {\bf B2-a}$(z_0, r_0)$ and {\bf B2-b}$(z_0, r_0)$ are satisfied for all $z_0\in \R^d$ and $r_0\in (0, \infty)$.
}
\end{example}

\begin{example}\label{ex:sdmms}
{\rm Suppose that $(\X, d, m)$ is an unbounded Ahlfors regular $n$-space for some $n>0$.
Assume that $d$ is uniformly equivalent to the shortest-path metric in $\X$. Suppose that
there is a diffusion process $Z$ with a symmetric, continuous transition density
$p^Z_t(x, y)$ satisfying the following sub-Gaussian bounds
\begin{align}\label{e:subgaussian}
&\frac{c_1}{t^{n/d_w}}\exp\left(-c_2\left(\frac{(d(x, y)^{d_w}}{t}\right)^{1/(d_w-1)}\right)
\le p^Z_t(x, y) \nonumber\\
&\,\,\,\, \le \frac{c_3}{t^{n/d_w}}\exp\left(-c_4\left(\frac{(d(x, y)^{d_w}}{t}\right)^{1/(d_w-1)}\right),
\end{align}
for all $x, y\in \X$ and $t>0$. Here $d_w\ge 2$ is the walk dimension of the space $\X$.
Let $T$ be a subordinator, independent of $Z$, with Laplace exponent $\phi$.
We define a process $X$ by $X_t=Z_{T_t}$. Then $X$ is a symmetric Feller process
and Assumption {\bf A} is clearly satisfied.

In this example, we will assume that the Laplace exponent $\phi$ is a complete Bernstein
function satisfying the following assumption:

\vspace{.1in}
\noindent
{\bf H}: There exist constants $0<\delta_1\le \delta_2 <1$ and $a_1, a_2>0$  such that
\begin{equation}\label{e:scaling4laplaceexp}
a_1\left(\frac{t}{s}\right)^{\delta_1} \le \frac{\phi( t)}{\phi( s)} \le a_2 \left(\frac{t}{s}\right)^{\delta_2} , \quad t\ge s > 0\, .
\end{equation}
Under this condition, by using \cite[Corollary 2.4 and Proposition 2.5]{KSV14a}
and repeating the argument of \cite[Lemmas 3.1--3.3]{KSV14a}
we obtain that that the jumping intensity $J$ of $X$ satisfies
$$
J(x, y)\asymp d(x, y)^{-n}\phi(d(x, y)^{-d_w}), \qquad x, y\in \X
$$
and that, when $2\delta_2<n/d_w$, $X$ is transient and its Green function $G$ satisfies
$$
G(x, y)\asymp d(x, y)^{-n}(\phi(d(x, y)^{-d_w}))^{-1}, \qquad x, y\in \X.
$$
Using these, one can easily see that Assumptions {\bf C},
and {\bf C1}$(z_0, r_0)$, {\bf C2}$(z_0, r_0)$,
{\bf D1}$(z_0, r_0)$ and {\bf D2}$(z_0, r_0)$  are satisfied for all $z_0\in \X$ and $r_0\in (0, \infty)$.
By using \cite[Proposition A.3]{BKuK} and repeating the argument of \cite[Corollary A.4]{BKuK},
one can easily show that Assumptions {\bf B}, {\bf B1-a}$(z_0, r_0)$,  {\bf B1-b}$(z_0, r_0)$,
{\bf B1-c}$(z_0, r_0)$, {\bf B2-a}$(z_0, r_0)$ and {\bf B2-b}$(z_0, r_0)$ are satisfied for all $z_0\in \X$ and $r_0\in (0, \infty)$.
}
\end{example}

\begin{example}\label{ex:sbm-like}
{\rm
Suppose that $T$ is a subordinator with Laplace exponent $\phi$.
In this example, we will assume that the Laplace exponent $\phi$ is a complete Bernstein
function satisfying the assumption {\bf H} in the previous example.

Suppose that $W$ is a Brownian motion in $\R^d$, independent of $T$, with generator $\Delta$.
The process $X$ defined by $X_t=W_{T_t}$ is a subordinate Brownian motion. $X$ is a symmetric
L\'evy process with L\'evy exponent $\phi(|\xi|^2)$ and its L\'evy density is given by $J_0(x)=
j_0(|x|)$ with
$$
j_0(r)=\int^\infty_0(4\pi t)^{-d/2}e^{-\frac{r^2}{4t}}\mu(t)dt, \qquad r>0,
$$
where $\mu(t)$ is the L\'evy density of $T$. It follows from \cite[Theorem 3.4]{KSV14a} that
\begin{equation}\label{e:2siededbd4densitysbm}
j_0(r)\asymp r^{-d}\phi(r^{-2}), \qquad r>0.
\end{equation}

Suppose that $k(x, y)$ is a symmetric function on $\R^d\times\R^d$ which is bounded between
two positive constants. The symmetric form $({\cal E}, C^2_c(\R^d))$ on $L^2(\R^d)$ defined by
$$
{\cal E}(f, g)=\int_{\R^d}\int_{\R^d}(f(x)-f(y))(g(x)-g(y))k(x, y)j_0(|x-y|)dxdy
$$
is closable, and so its minimal extension $({\cal E}, {\cal F})$ is a regular Dirichlet
form. Let $X$ be the symmetric Markov process on $\R^d$ associated with $({\cal E}, {\cal F})$.
It follows from \cite{CK08} that $X$ is a conservative Feller process and it admits a
transition density $p(t, x, y)$ satisfying the following estimates: for all $(t, x, y)\in
(0, \infty)\times\R^d\times\R^d$,
\begin{equation}\label{e:sbmlikedensity}
c^{-1}\left(\Phi^{-1}(t)^d\wedge (tj_0(|x-y|))\right)\le p(t, x, y)
\le c\left(\Phi^{-1}(t)^d\wedge (tj_0(|x-y|))\right),
\end{equation}
where $\Phi^{-1}$ is the inverse of the function
$$
\Phi(r)=\frac1{\phi(r^{-2})}, \qquad r>0.
$$
When $d>2\delta_2$, $X$ is transient and its Green function $G(x, y)$ satisfies
$$
G(x, y)\asymp \frac1{|x-y|^d\phi(|x-y|^{-2})}, \qquad x\neq y\in \R^d.
$$

If we further assume that the first partial derivatives of $k$ are bounded and continuous on
$\R^d\times \R^d$, then for $f\in C^2_c(\R^d)$, the generator ${\cal A}$ of $X$
admits the following expression
\begin{align}\label{e:snmlikegen}
{\cal A}f(x)=&\int_{\R^d}(f(x+y)-f(x)-y\cdot\nabla f(x)1_{|y|<1})k(x, x)j_0(|y|)dy\nonumber\\
&+\int_{\R^d}(f(x+y)-f(x))(k(x, x+y)-k(x, x))j_0(|y|)dy.
\end{align}

Using \eqref{e:2siededbd4densitysbm}--\eqref{e:snmlikegen}, one can easily check that
Assumptions {\bf A}
{\bf B}, {\bf B1-a}$(z_0, r_0)$,  {\bf B1-b}$(z_0, r_0)$,
{\bf B1-c}$(z_0, r_0)$, {\bf B2-a}$(z_0, r_0)$, {\bf B2-b}$(z_0, r_0)$,
{\bf C}, {\bf C1}$(z_0, r_0)$, {\bf C2}$(z_0, r_0)$,  {\bf D1}$(z_0, r_0)$ and {\bf D2}$(z_0, r_0)$
are satisfied for all $z_0\in \R^d$ and $r_0\in (0, \infty)$.
}
\end{example}

\section{Limit of Green functions at regular boundary points}\label{s:green}

Suppose that $D\subset \X$ is an open set.
In this section, we will prove that, under some assumptions,
the Green function $G_D(x, y)$ of
$X^D$ approaches zero when $x$ approaches a point $z\in \partial D$
which is regular for $D^c$ with respect to $X$.

\begin{prop}\label{p:cont_b_1}
Suppose that
$X$ is a Hunt process on $\X$ satisfying both the Feller and
the strong Feller property.
If $z$ is a regular boundary point
of $D$ and $f$ is a bounded
Borel function on $ D^c$ which is continuous at $z$, then
$$
\lim_{\overline{D}\ni x\to z}\E_x\left[f(X_{\tau_D});
\tau_D<\infty \right]=f(z).
$$
\end{prop}

\pf
Note that, since  $z$ is a regular boundary point, $\P_z(\tau_D=0)=1$.
By \cite[Lemma 3]{chung}, for every $s>0$,
$$
\limsup_{x\to z} \P_x(\tau_D>s)  \le \P_z(\tau_D>s)  =0.
$$
Thus, for every $s>0$,
\begin{equation}\label{e:pcont_b}
\lim_{x\to z} \P_x(\tau_D>s) =0.
\end{equation}
In particular,
$\lim_{x\to z} \P_x(\tau_D<\infty) =1$.
Thus it is enough to show that,
\begin{equation}\label{e:pcont_bpp}
 \lim_{\overline{D}\ni x\to z}\E_x\left[|f(X_{\tau_D})-f(z)| :
\tau_D<\infty \right] < \eps
\end{equation}
for arbitrary $\eps>0$.

Given $\eps>0$, choose $\delta>0$ such that
$$|f(w)-f(z)| \le \frac{\eps}{2}, \quad \text {for every } w \in B(z, \delta).$$
Now we have for every $s>0$ and $x \in \overline{B(z, \delta/2)}$,
\begin{eqnarray*}
\lefteqn{\E_x\left[|f(X_{\tau_D})-f(z)| :
\tau_D<\infty \right]}\\
 & =&\E_x\left[|f(X_{\tau_D})-f(z)| :
\tau_D<\tau_{B(z, \delta)} \right] + \E_x\left[|f(X_{\tau_D})-f(z)| :
\tau_{B(z, \delta)} \le \tau_D<\infty \right] \\
&\le &\frac{\eps}{2} \P_x(
\tau_D<\tau_{B(z, \delta)} ) +2 \|f\|_\infty  \P_x(
\tau_{B(z, \delta)} \le s \text{ or } s <\tau_D ) \\
&\le &\frac{\eps}{2} +2 \|f\|_\infty  (\P_x(
\tau_{B(z, \delta)} \le s )+\P_x(
 s <\tau_D ) ).
\end{eqnarray*}
It follows from \cite[Lemma 2]{chung} that there exists an $s>0$ such that
$$
\sup_{x \in \overline{B(z, \delta/2)}}\P_x(
\tau_{B(z, \delta/2)} \le s )\le \frac{\eps}{4 \|f\|_\infty}.
$$
Using \eqref{e:pcont_b}, we get
$$
\lim_{\overline{D}\ni x\to z}\E_x\left[|f(X_{\tau_D})-f(z)| :
\tau_D<\infty \right] < \eps +2 \|f\|_\infty \lim_{\overline{D}\ni x\to z}\P_x(
 s <\tau_D ) =\eps .
 $$
 We have proved \eqref{e:pcont_bpp}.
\qed

The following result will be used in \cite{KSV15a}.

\begin{prop}\label{p:contgen00}
Suppose that $X$ is a Hunt process on $\X$ satisfying both the Feller and
the strong Feller property and that for all Greenian open sets $V$,
$x\mapsto G_{V} (x,y)$ is continuous in $V \setminus\{y\}$.
If $D\subset\X$ is an open set, $z_1\in \partial D$ is regular for $D^c$
and there exists $r_0>0$ such that $D\cup B(z_1, r_0)$ is Greenian, then for all $y \in D$,
$$
\lim_{D\ni x\to z_1}G_D(x, y)=0.
$$
\end{prop}

\pf
Fix $y \in D$ and choose $r_1 \le r_0/2$ small enough so that $y \in B(z_1, 4 r_1)^c \cap D$.
Let $U_1=D\cup B(z_1, r_1)$
and $U_2=D\cup B(z_1, 2r_1)$
which are both Greenian. Then, by our assumption, $x\mapsto G_{U_i} (x, y)$ are
continuous in $U_{i}\setminus\{y\}$.

By the strong Markov property we have
$$
G_D(x,y) = G_{U_1}(x,y) - \E_x [ G_{U_1}(X_{\tau_D} , y)].
$$
The function $w\mapsto G_{U_1} (w, y)$  is bounded on $D^c$.
Indeed, using domain-monotonicity of Green functions and continuity of  $G_{U_2}(\cdot, y)$ on  $U_{2}\setminus\{y\}$,
$$
\sup_{w \in D^c} G_{U_1}(w, y)  \le
\sup_{w \in B(z_1, r_1)} G_{U_1}(w, y)  \le \sup_{w \in B(z_1, r_1)} G_{U_2}(w, y) <\infty.
$$
Since $x\mapsto G_{U_1}(x, y)$ is continuous at $z_1$, by Proposition \ref{p:cont_b_1},
we have
$$
\lim_{D\ni x\to z_1}G_D(x, y)= G_{U_1}(z,y)- \lim_{D\ni x\to z_1}\E_x [ G_{U_1}(X_{\tau_D} , y)]    =0.
$$
\qed

The following result is quite general.

\begin{prop}\label{p:contharm}
Suppose that $X$ is a Hunt process on $\X$ satisfying both the Feller and
the strong Feller property.

\noindent
(a) If $U\subset \X$ is open and $u$ is a bounded Borel
function on $U^c$, then the function $x\mapsto\E_x[u(X_{\tau_U}); \tau_U<\infty]$ is continuous in $U$.

\noindent
(b)
Assume further
that $X$ satisfies the Harnack
principle in the sense that, for any $z_1$,
there exists $r_0>0$ and $c>0$ such that for every $r\in (0, r_0)$ and every function
$u$ which is nonnegative on $\X$ and harmonic in $B(z_1, r)$ with respect to $X$,
it holds that
$$
u(x)\le cu(y), \qquad x, y\in B(z_1, r/2).
$$
Then, if $h$ is a nonnegative function on $\X$ which is harmonic in an open set $D\subset\X$ with respect to $X$,
then $h$ is continuous in $D$.
\end{prop}

\pf Part (a) follows from \cite[Theorem 3.4]{sch-wan12}.
Note that although \cite[Theorem 3.4]{sch-wan12} is stated for the case $\X=\R^d$, the argument there works generally.
Now one can repeat the argument after \cite[Theorem 2.1]{KSV15} to get the conclusion of (b).
\qed

\begin{corollary}\label{c:contgen00}
Suppose that $X$ is a Hunt process on $\X$ satisfying both the Feller and
the strong Feller property and  $X$ satisfies the Harnack
principle.
 Assume that $D\subset\X$ is an open set,
 that $z_1\in \partial D$ is regular for $D^c$
and there exists $r_0>0$ such that $D\cup B(z_1, r_0)$ is Greenian. Then for all $y \in D$,
$$
\lim_{D\ni x\to z_1}G_D(x, y)=0.
$$
\end{corollary}

We now weaken the Greenian assumption.

\begin{prop}\label{p:contgen0}
Suppose that $X$ satisfies Assumption {\bf A} and that, for every $z_1\in \X$, there
is $r_0>0$ such that the conclusion of Corollary \ref{c:bhp-finite} (BHP)
holds.
Assume that $D\subset\X$ is an open Greenian set,
 that $z_1\in \partial D$ is regular for $D^c$
and the
open balls $B(z_1, r)$
 are Greenian for all $r>0$. Then for all $y \in D$,
$$
\lim_{D\ni x\to z_1}G_D(x, y)=0.
$$
\end{prop}

\pf
Fix $y \in D$ and let $r_1=2d(z_1, y)$ and
 $U=D\cap B(z_1, r_1)$ which is Greenian.
By the strong Markov property we have
$$
G_D(x,y) = G_U(x,y) + \E_x [ G_D(X_{\tau_U} , y)].
$$
The BHP assumption implies that $X$ satisfies the Harnack
principle. Moreover, the open set  $U \cup B(z_1, r_1)=B(z_1, r_1)$  is Greenian.
Therefore we can apply Corollary \ref{c:contgen00} and get
$$
\lim_{D\ni x\to z_1}G_U(x, y)= \lim_{U\ni x\to z_1}G_U(x, y)  =0.
$$
Define
$$
w(x)=\P_x[X_{\tau_U}\in U^c\setminus B(z_1, r_1); \tau_U<\infty], \qquad x\in U.
$$
It follows from Proposition \ref{p:cont_b_1} that
\begin{align}\label{e:fins}
\lim_{U\ni x\to z_1}w(x)=0.
\end{align}
Note that $x  \to \E_x [ G_D(X_{\tau_U} , y)]$ and $w$ are both regular harmonic in $U$, vanish in $B(z_1, r_1)\cap(\overline{U}^c\cup
U^{{\mathrm reg}})$. Now we can combine \eqref{e:fins} with the BHP to get
$$
\lim_{U\ni x\to z_1}\E_x [ G_D(X_{\tau_U} , y)] =0.
$$
Therefore
$$
\lim_{D\ni x\to z_1}G_D(x,y) = \lim_{U\ni x\to z_1}G_U(x,y) + \lim_{U\ni x\to z_1}\E_x [ G_D(X_{\tau_U} , y)] =0.
$$
\qed

\section{Appendix}\label{s:appendix}
In this section we present the proof of Theorem \ref{t:main-finite}.
Throughout this section, $z_0$ is a fixed point in $\X$.
We will always assume in this section that, in addition to {\bf A}, {\bf B} and {\bf C}, the assumptions
\eqref{e:V-doubling}, \eqref{e:V-another1},  \eqref{e:Phi-doubling}, {\bf B1-a}$(z_0, r_0)$, {\bf B1-b}$(z_0, r_0)$, {\bf B1-c}$(z_0, r_0)$, {\bf C1}$(z_0, r_0)$ and {\bf D1}$(z_0, r_0)$
hold true for some $r_0\in (0,  R_0]$.
Recall that $n_0$ is the natural number in \eqref{e:V-another1}
.
In the next result we understand $r_0/(2n_0)$ to be $\infty$ if $r_0=\infty$.
\begin{prop}\label{p:pcom-fin}
There exists a constant $c=c(z_0, r_0)>0$ such that for all
$r<r_0/(2n_0)$ and all $x\in B(z_0,r)$,
$\E_x \tau_{B(z_0,r)} \le c \Phi(r)$.
\end{prop}
\pf
Let
$r<r_0/(2n_0)$, denote $B=B(z_0,r)$ and let $F(t):=\P_x(\tau_B>t)$.
First note that if $y\in B$, then by {\bf C1}$(z_0, r_0)$, $j(y,z)\ge c_1 j(z_0,z)$ for all $z\in \overline{A}(z_0,2r,
2n_0 r)$. Hence,
$$
J(y,\X \setminus B)\ge J(y, \overline{A}(z_0,2r,
2n_0r))\ge c_1 J(z_0, \overline{A}(z_0,2r,
2n_0r))\, .
$$
In the same way as in \cite[Proposition 2.1]{BKuK}, $-F'(t)\ge c_1  J(z_0, \overline{A}(z_0,2r,
2n_0r))$, implying that $F(t)\le \exp\{-t c_1  J(z_0, \overline{A}(z_0,2r,
2n_0r)\}$. Hence,
\begin{equation}\label{e:pcom-fin}
\E_x \tau_{B(z_0,r)} \le \left(c_1  J(z_0, \overline{A}(z_0,2r,
2n_0r))\right)^{-1}\, .
\end{equation}
By using {\bf C1}$(z_0, r_0)$ and the monotonicity of $V$ and $\Phi$ in the first line, \eqref{e:V-doubling} and \eqref{e:Phi-doubling} in the second line, \eqref{e:V-another1}
  in the third, we get
\begin{eqnarray*}
J(z_0,A(z_0,2r,
2n_0r))&=& \int_{\overline{A}(z_0,2r,2n_0r)} j(z_0,z)m(dz) \ge  \int_{\overline{A}(z_0,2r,2n_0r)} \frac{c_2}{V(n_0r)\Phi(n_0r)} m(dz)\\
&\ge & \int_{\overline{A}(z_0,2r,2n_0r)} \frac{c_3}{V(2r)\Phi(r)} m(dz) \ge \frac{c_3}{\Phi(r)}\left(\frac{V(2n_0r)}{V(2r)}-1\right)\\
&\ge &\frac{c_4}{\Phi(r)}.
\end{eqnarray*}
Together with \eqref{e:pcom-fin} this proves the claim. \qed

Let $a\in (1/2,1)$. For each $r <r_0$, we consider a function $\ph^{(r)}\in \dom(\overline{B}(z_0, ar), B(z_0, r))$,
and let  $V^{(r)}=\{x\in \X: \ph^{(r)}(x)>0\}$.
Note that, by choosing $\ph^{(r)}$ appropriately, we can achieve that
$\delta^{(r)}:= \sup_{x \in \X} \max\big(\sA \phi^{(r)}(x), \widehat{\sA} \phi^{(r)}(x)\big)\le {c}/{\Phi(r)}$,
where $c=c(z_0, r_0, a)$ is the constant in assumption {\bf B1-a}$(z_0, r_0)$.

In what follows, our analysis and results are valid for all $r< r_0$ with constants depending on $a\in (1,2)$,
but \emph{not} on $r$. To ease the notation in the remaining
part of the section we drop the superscript $r$ from $\ph^{(r)}$ and $V^{(r)}$ and write simply $\ph$ and $V$.

Let
\begin{align}\label{e:defpsi-fin}
 \psi(x)  = \frac{\max(\sA \ph(x), \widehat{\sA} \ph(x), \delta (1 - \ph(x)))}{\ph(x)} \, ,
  \quad   x \in \X ,
\end{align}
with the convention $1/0=\infty$. Note that $\psi(x)=\infty$ for $x\in V^c$, and $\psi(x)=0$ for $x\in \overline{B}(z_0,ar)$.

As in Section \ref{s:proofs}
we define two right-continuous additive functionals $A_t$ and $\widehat A_t$
as in \eqref{e:eAt}
and define two right-continuous exact  strong Markov multiplicative functionals
$M_t = \exp(-A_t)$ and $\widehat M_t = \exp(-\widehat A_t)$. We consider the semigroup of operators
$T^\psi_t f(x) = \E_x(f(X_t) M_t)$
associated with the multiplicative functional $M$, which is the transition operator of subprocess $X^\psi$ of $X$
and the semigroup of operators $\widehat T^\psi_t f(x) = \E_x(f(\widehat X_t) \hat M_t)$
associated with the multiplicative functional $\widehat M$, which is the transition operator of subprocess $\widehat{X}^\psi$ of $\widehat{X}$.
Again, the  potential densities of $X^\psi$ and $\widehat{X}^\psi$
satisfy $\widehat{G}^\psi(x, y) = G^\psi(y, x)$
and
\begin{align}\label{e:Gmon-fin}
 G^\psi(x,y) \le  G_V(x,y)  \le  G_{B(z_0, r)} (x,y),  \quad (x,y) \in V \times V.
\end{align}

Let $\tau_a = \inf \set{t \ge 0 : A_t \ge a}$ and
\begin{equation}\label{e:psi-tau-a-fin}
 \pi^\psi f(x) =-\E_x \int_{[0, \infty)} f(X_t) d M_t = \E_x \expr{\int_0^\infty f(X_{\tau_a}) e^{-a} da} = \int_0^\infty \E_x(f(X_{\tau_a})) e^{-a} da.
\end{equation}
Recall from \cite[pp.~492--493]{BKuK} that $\pi^\psi f$ can be written in the following two ways:
if $f$ is nonnegative and vanishes in $\X \setminus (\overline{B}(z_0, ar)^c \cap V)$, then
\begin{align}\label{e:pipsi1-fin}
 \pi^\psi f(x) & = G^\psi (\psi f)(x) = \int_{V \cap \overline{B}(z_0, ar)^c} G^\psi(x, y)
 \psi(y) f(y) m(dy), \quad  x \in   V,
\end{align}
and if $f \in \dom(\sA)$ vanishes in $V$ then for all $x \in V$,
\begin{align}\label{e:pipsi2-fin}
& \pi^\psi f(x)  = G^\psi \sA f (x) = \int_V G^\psi(x, y) \sA f(y) m(dy)\nonumber\\
&=\int_V G^\psi(x, y) \int_{\X \setminus V} f(z) j(y, z)  m(dz) m(dy)\nonumber\\
 &=\int_{\X \setminus V} \expr{\int_V G^\psi(x, y) j(y, z) m(dy)} f(z) m(dz).
\end{align}
By Corollary 4.8 of \cite{BKuK}, we have
\begin{align} \label{e:Vzero-fin}
\pi^\psi(x, \partial V) = 0, \qquad x\in V.
\end{align}
By \cite[Lemma 4.10]{BKuK}, we get that if $f$ is regular harmonic
in $B(z_0, r)$ with respect to $X$,
then  $f(x) = \pi^\psi f(x)$
for all $x\in \overline{B}(z_0, r/2)$.
The main step of the proof is to get the correct estimate of $\pi^\psi(x, dy)/m(dy)$.

Let $U$ be an open subset of $V$. For any nonnegative or bounded $f$ and $x \in V$ we let
\begin{align*}
 \pi^\psi_U f(x) & = \E_x (f(X_{\tau_U}) M_{\tau_U-}) , & G^\psi_U f(x) & = \E_x \int_0^{\tau_U} f(X_t) M_t dt .
\end{align*}
 $G^\psi_U$ admits a density $G^\psi_U(x, y)$, and we have $G^\psi_U(x, y) \le G_U(x, y)$, $G^\psi_U(x, y) \le G^\psi(x, y)$.
 For any $f \in \dom(\sA)$ we have
\begin{align}
\label{e:psid2-fin}
 \pi^\psi_U f(x) & = G^\psi_U (\sA - \psi) f(x) + f(x) , && x \in V .
\end{align}
In particular, by an approximation argument,
\begin{align}
\label{e:psilevy-fin}
  \pi^\psi_U(x, E) & = \int_U G^\psi_U(x, y) J(y, E) m(dy) , && x \in U , \, E \sub \X \setminus \overline{U} .
\end{align}

By the definition of $\psi$, we have that $(\sA - \psi) \ph(x) \le 0$ for $x \in \X$ (and, in particular, for all $x\in V$). Thus using \eqref{e:psid2},
we have the following.

\begin{lemma}[{\cite[Lemma 4.4]{BKuK}}]
\label{l:4.4-fin}
Let $U = V \cap \overline{B}(z_0, ar)^c$. Then
\begin{align}
\label{e:4.4-fin}
 \pi^\psi_U(x, V\setminus U) & \le \ph(x) , && x \in U .
\end{align}
\end{lemma}

\begin{lemma}\label{l:4.5-fin}
Let $b\in (1/2,a)$. There exists a constant  $c=c(a,b)>0$ such that
\begin{align*}
 G^\psi(x, y) & \le c  \frac{\Phi(r)}{V(r)} \ph(x), \quad   x \in V \cap B(z_0, br)^c , \, y \in B(z_0,r/2)\, .
\end{align*}
\end{lemma}
\pf If $x\in \overline{A}(z_0,br,ar)$, then $\ph(x)=1$. Thus, by Assumption {\bf D1}$(z_0, r_0)$,  for $x \in \overline{A}(z_0, br, ar)\subset \overline{A}(z_0, br,r)$ and $y \in B(z_0,r/2)$,
\begin{align*}
 G^\psi(x, y)  \le  G_V(x, y) \le  G_{B(z_0, r)} (x, y)   \le c_1  \frac{\Phi(r)}{V(r)}= c_1  \frac{\Phi(r)}{V(r)} \ph(x),
\end{align*}
with $c_1=c_1(b)$.

For the remainder of the proof we assume that $U:=V\cap \overline{B}(z_0,ar)^c$.  Let $f\ge 0$ be supported in $\overline{B}(z_0,r/2)$ with $\int f(w)m(dw)=1$. Then, by the strong Markov property,
$$
G^\psi f(x) =\pi^\psi_U (G^\psi f)(x)=\pi^\psi_U ( {\bf 1}_{\overline{A}(z_0, br, ar)}G^\psi f)(x)+\pi^\psi_U ( {\bf 1}_{B(z_0, br)} G^\psi f)(x)=:I+II.
$$
First note that by  {\bf D1}$(z_0, r_0)$ and \eqref{e:Gmon-fin}, for $y \in \overline{A}(z_0, br, ar)$,
$$
G^\psi f(y) \le \int_{\overline{B}(z_0, r/2)} G_V(y,w) f(w)m(dw)  \le \int_{\overline{B}(z_0, r/2)} G_{B(z_0, r)}(y,w) f(w)m(dw) \le c_2 \frac{\Phi(r)}{V(r)}.
$$
Thus, combining  this with Lemma \ref{l:4.4-fin} we get
\begin{align*}
&I \le \left(\sup_{y \in \overline{A}(z_0, br, ar)} G_\psi f(y) \right)
  \pi^\psi_U(x, \overline{A}(z_0, br, ar)) \le   c_2\frac{\Phi(r)}{V(r)}  \pi^\psi_U(x, B(z_0,ar))  \le c_2  \frac{\Phi(r)}{V(r)} \ph(x).
\end{align*}
For $II$, note that by {\bf C1}$(z_0, r_0)$, for every $z\in B(z_0,br)$,
$$
\int_U G_U^{\psi}(x,y)j(y,z_0)\, m(dy)\le c_3 \int_U G_U^{\psi}(x,y)j(y,z)\, m(dy)
$$
for a constant $c_3=c_3(a,b)$. By integrating over the ball $B(z_0,br)$, using the doubling property of $V$,
Lemma \ref{l:4.4-fin} and \eqref{e:psilevy-fin},  we obtain that
\begin{eqnarray*}
\int_U G_U^{\psi}(x,y)j(y,z_0)m(dy) &\le & \frac{c_4}{V(r)}\int_{B(z_0,br)} \int_U G_U^{\psi}(x,y)j(y,z)\, m(dy) \, m(dz)\\
&= & \frac{c_4}{V(r)} \int_U G_U^{\psi}(x,y)\left(\int_{B(z_0,br)}j(y,z)\, m(dz)\right) m(dy)\\
&=& \frac{c_4}{V(r)} \pi_U^{\psi}(x,B(z_0,br)) \le \frac{c_4}{V(r)} \ph(x)\,  ,
\end{eqnarray*}
with $c_4=c_4(a,b)$. Thus, by using {\bf C1}$(z_0, r_0)$ in the third line and the display above in the last line,
\begin{eqnarray*}
II &= & \int_U G_U^{\psi}(x,y) \int_{B(z_0,br)} j(y,z) G^{\psi}f(z)\, m(dz)\, m(dy) \\
&=& \int_{B(z_0,br)}G^{\psi}f(z)\left( \int_U G_U^{\psi}(x,y)j(y,z)\, m(dy)\right) m(dz)\\
&\le & c_5 \int_{B(z_0,br)}G^{\psi}f(z)\left( \int_U G_U^{\psi}(x,y)j(y,z_0)\, m(dy)\right) m(dz)\\
&\le & \frac{c_6}{V(r)}\ph(x) \int_{B(z_0,br)}G^{\psi}f(z)\, m(dz)\,  .
\end{eqnarray*}
Finally, by using the dual version of Proposition \ref{p:pcom-fin},
\begin{eqnarray*}
\int_{B(z_0,br)}G^{\psi}f(z)\, m(dz)&\le & \int_{B(z_0,r/2)} \left(\int_{B(z_0,r)} G_V(z,y)\, m(dz)\right) f(y)\, m(dy)\\
&\le & \int_{B(z_0,r/2)} \left(\int_{B(z_0,r)} G_{B(z_0,r)}(z,y)\, m(dz)\right) f(y)\, m(dy)\\
&=& \int_{B(z_0,r/2)} \widehat{\E}_y (\tau_{B(z_0,r)}) f(y)\, m(dy) \le c_7 \Phi(r)\,  .
\end{eqnarray*}
This completes the proof. \qed

\begin{lemma}\label{l:4.9-fin}
There exists a constant $c=c(a,z_0)>0$ such that for all $x \in B(z_0,r/2)$,
\begin{align}
\label{e:pipsiest-fin}
\pi^\psi(x, dy)/m(dy) \le  c\Phi(r) j(z_0,y)\ind_{\overline{B}(z_0, ar)^c}(y)\, .
\end{align}
\end{lemma}
\pf
Let $b:=\frac{2a}{1+2a}$ so that $b\in (1/2,a)$. First note that $\psi$ vanishes on $\X \setminus \overline{B}(z_0, ar)$.
Thus
\begin{align} \label{e:Vsupp-fin}
\pi^\psi(y, \overline{B}(z_0, ar)) = 0, \quad y \in V.
\end{align}
Fix $x \in B(z_0,r/2)$. If $f$ is a non-negative function on $\X$ vanishing in $\X \setminus (\overline{B}(z_0, ar)^c \cap V)$, then by \eqref{e:pipsi1-fin} and the dual version of Lemma \ref{l:4.5-fin} (together with $G^{\psi}(x,y)=\widehat{G}^{\psi}(y,x)$),
\begin{align}\label{e:pipsi11-fin}
 \pi^\psi f(x)  \le   c_1  \frac{\Phi(r)}{V(r)} \int_{V \cap \overline{B}(z_0, ar)^c} \ph(y)    \psi(y) f(y) m(dy).
\end{align}
Since for $y\in \overline{B}(z_0,ar)^c$ we have $\ph(y)    \psi(y)  \le c_2(\Phi(r))^{-1}$ by the definition in
\eqref{e:defpsi-fin}, assumption {\bf B1-a}$(z_0, r_0)$
and \eqref{e:AC21-finite} (note that $V\cap \overline{B}(z_0, ar)^c\subset \overline{A}(z_0,ar,r)$), we get
\begin{align}\label{e:pipsi12-fin}
\pi^\psi f(x)  \le    \frac{c_3}{V(r)}\int_{V \cap \overline{B}(z_0, ar)^c} f(y) m(dy) \le c_4 \Phi(r)\int_{V \cap \overline{B}(z_0, ar)^c} j(z_0,y)f(y)\, m(dy)\, .
\end{align}
On the other hand, if $g \in \dom(\sA)$ vanishes in $V$ then by \eqref{e:pipsi2-fin},
\begin{align}\label{e:pipsi21-fin}
\pi^\psi g(x)   = \int_{\X \setminus V} \expr{\int_V G^\psi(x, y) j(y, z) m(dy)} g(z) m(dz).
\end{align}
Assume $z \in \X \setminus V \subset \overline{B}(z_0, ar)^c$ and let
$$
I:= \int_{V \cap  \overline{B}(z_0, br)^c}  G^\psi(x, y) j(y, z) m(dy)\quad \text{and}
\quad II:=\int_{B(z_0, br)} G^\psi(x, y) j(y, z) m(dy).
$$
We now consider $I$ and $II$ separately.

By the dual versions of Lemma \ref{l:4.5-fin} and \eqref{e:Assumptions B1(b)}, and the fact that $\widehat{\sA}\ph(z)=\widehat{J}\ph(z)$, for  $c_4=c_4(a)>0$
\begin{align}\label{e:pipsi23-fin}
I \le c_4  \frac{\Phi(r)}{V(r)}   \int_{V \cap  \overline{B}(z_0, br)^c}  \ph(y)    j(y, z) m(dy)
= c_4  \frac{\Phi(r)}{V(r)}  \widehat  \sA \ph(z) \le  c j(z_0, z) \Phi(r)\, .
\end{align}

On the other hand, by assumption {\bf C1}$(z_0, r_0)$ and \eqref{e:Gmon-fin}, for $c_6=c_6(a)>0$
\begin{align}\label{e:pipsi24-fin}
II &\le c_5 \int_{B(z_0, br)} G^\psi(x, y) j(z_0, z) m(dy)\nonumber\\
 &\le c_5 j(z_0,z)\int_{B(z_0, r)} G_{B(z_0,r)}(x, y) m(dy)\nonumber\\
 &=  c_5 j(z_0,z) \E_x \tau_{B(z_0,r)}\le c_6 j(z_0,z) \Phi(r)\,  .
\end{align}
Hence,
$$
\pi^{\psi} g(x)\le c_7 \Phi(r)\int _{\X\setminus V}j(z_0,z)g(z)\, m(dz)\, .
$$
Together with \eqref{e:pipsi12-fin} this proves the lemma. \qed

\begin{corollary}\label{c:4.9-fin}
Let $f$ be a non-negative function on $\X$ and $x\in B(z_0,r/2)$ such that $f(x)\le \E_x f(X_{\tau})$ for every stopping time $\tau\le \tau_{B(z_0, r)}$. Then
\begin{equation}\label{e:har:c-fin}
f(x) \le c\Phi(r) \int_{\overline{B}(z_0, ar)^c} j(z_0,y)f(y) m(dy)\,  ,
\end{equation}
where $c=c(a)$ is the constant from Lemma \ref{l:4.9-fin}.
\end{corollary}
\pf
Recall from \eqref{e:psi-tau-a-fin} that $\pi^\psi f(x)=\int_0^\infty \E_x(f(X_{\tau_a})) e^{-a} da$. Since $\tau_a\le \tau_V\le \tau_{B(z_0,r)}$, we have that $f(x) \le \E_x f(X_{\tau_a})$, and therefore $f(x) \le \pi_\psi f(x)$. Thus by \eqref{e:pipsiest-fin},
$$
f(x) \le \int f(y) \pi^\psi(x, dy) \le c\Phi(r) \int_{\overline{B}(z_0, ar)^c} j(z_0,y) f(y) m(dy)\, .
$$
\qed

\begin{lemma}\label{l:exit-fin}
For any $b\in (1/2,a)$ there exists $c=c(z_0, r_0, a,b)>0$ such that for any $r < r_0$ and any open set
$D\subset B(z_0, br)$ we have
$$
{\P}_x\left(X_{\tau_D} \in \overline{A}(z_0, br,ar)\right) \,\le\,  \frac{c}{\Phi(r)} \E_x\tau_D\, , \qquad x \in D\cap B(z_0, r/2)\, .
$$
\end{lemma}
\pf Let $f\in \dom(\overline{A}(z_0, br, ar), A(z_0, r/2, r))$. By assumption {\bf B1-c}$(z_0, r_0)$, $\sup_{y\in\X}\sA f(y)\le \frac{c}{\Phi(r)}$ with $c=c(z_0,a,b)$. By Dynkin's formula, for $x\in D\cap B(z_0,r/2)$,
$$
\E_x f(X_{\tau_D})=\E_x\int_0^{\tau_D} \sA f(X_t)dt \le \frac{c}{\Phi(r)} \E_x \tau_D\, .
$$
The claim follows from $\ind_{\overline{A}(z_0, br, ar)}\le f$. \qed

\begin{prop}\label{p:main-fin-aux}
Let $b\in (1/2,a)$. There exists $c=c(z_0, r_0,a,b)>1$ such that for any $r <r_0$, any
open set $D \subset B(z_0,r)$ and any non-negative function $u$ on $\X$ which is
regular harmonic with respect to $X$ in $D$ and
vanishes on $B(z_0,r)\cap \left(\overline{D}^c\cup D^{\mathrm{reg}} \right)$, it holds that
\begin{align}\label{e:bhp-inf-lemma-fin}
&c^{-1}\E_x \tau_{D\cap B(z_0,br/2)} \int_{\overline{B}(z_0, ar/2)^c} j(z_0,z)u(z)\, m(dz) \le u(x)\nonumber\\
 &\,\,\,\, \le \E_x \tau_{D\cap B(z_0,br/2)} \int_{\overline{B}(z_0, ar/2)^c} j(z_0,z)u(z)\, m(dz)
\end{align}
for all $x\in D\cap B(z_0,r/4)$.
\end{prop}
\pf
Let $O:=D \cap B(z_0, br/2)$, $D_1:=\overline{A}(z_0, br/2, ar/2)$ and $D_2:=B(z_0, ar/2)^c$.
By the harmonicity of $u$,
\begin{align}\label{e:bhp1-fin}
 u(x)=\E_x[u(X_{\tau_{O}})]=\E_x[u(X_{\tau_{O}}):X_{\tau_{O}}\in D_1 ]+\E_x[u(X_{\tau_{O}}):X_{\tau_{O}}\in D_2 ]\, , \ \ x\in D.
\end{align}
Since $u$ vanishes on $B(z_0,r)\cap \left(\overline{D}^c\cup D^{\mathrm{reg}} \right)$,
it follows that $u(y)\le \E_y u(X_{\tau})$ for every stopping time $\tau\le \tau_{B(z_0,r)}$
and every $y\in B(z_0, r)\setminus D^{\mathrm{irr}}$.
Since $D^{\mathrm{irr}}$ is polar with respect to $X$, we see that $X_{\tau_O}\notin D^{\mathrm{irr}}$.
It  follows from Corollary \ref{c:4.9-fin}
and Lemma \ref{l:exit-fin} that for all $x\in D\cap B(z_0,r/4)$,
\begin{align}\label{e:bhp2-fin}
&\E_x[u(X_{\tau_{O}}):X_{\tau_{O}}\in D_1 ]
\le \left(\sup_{y \in D_1\setminus D^{\mathrm{irr}}} u(y) \right) \P_x(X_{\tau_{O}}\in D_1) \nonumber\\
&\le   \frac{c_1}{\Phi(r)} (\E_x\tau_O) \Phi(r/2)\int_{\overline{B}(z_0,a^2r)^c}j(z_0,z)u(z)\, m(dz) \\
&\le c_1 \E_x\tau_O \int_{\overline{B}(z_0,ar/2)^c}j(z_0,z)u(z)\, m(dz)
\end{align}
with $c_1=c_1(a,b)$.
On the other hand, by assumption {\bf C1}$(z_0, r_0)$, for all $x\in D\cap B(z_0,r/4)$,
\begin{align}\label{e:bhp3-fin}
\E_x[u(X_{\tau_{O}}):X_{\tau_{O}}\in D_2]
&=\int_{\overline{B}(z_0, ar/2)^c} \int_{O} G_{O}(x,y)j(y,z)m(dy)u(z)m(dz)\nonumber\\
&\asymp
\int_{\overline{B}(z_0, ar/2)^c} \int_{O} G_{O}(x,y)j(z_0,z)m(dy)u(z)m(dz)\nonumber\\
&=\E_x \tau_O\int_{\overline{B}(z_0, ar/2)^c}j(z_0,z)u(z)m(dz).
\end{align}
The proposition now follows from \eqref{e:bhp1-fin}--\eqref{e:bhp3-fin}.
\qed

\begin{lemma}\label{l:exit-comparison}
For any $b\in (1/2,a)$ there exists $c=c(a,b)>0$ such that for every
$r < r_0/(4n_0)$, and every open set $D\subset B(z_0,2r)$,
$$
\E_x \tau_{D\cap B(z_0,br)}\le \E_x \tau_D \le c \E_x \tau_{D\cap B(z_0,br)}\,  , \qquad x\in D\cap B(z_0,ab r)\, .
$$
\end{lemma}
\pf First note that by the strong Markov property,
$$
\E_x \tau_D = \E_x \tau_{D\cap B(z_0,br)} +\E_x\big[\E_{X_{\tau_{D\cap B(z_0,br)}} }\tau_D \big]\,  .
$$
By Proposition \ref{p:pcom-fin}, Lemma \ref{l:exit-fin} and doubling property of $\Phi$, for $x\in D\cap B(z_0, ab r)$,
\begin{eqnarray*}
\E_x\big[\E_{X_{\tau_{D\cap B(z_0,br)}} }\tau_D \big] &\le & \big(\sup_{y\in D} \E_y \tau_D\big)  \P_x \big(X_{\tau_{D\cap B(z_0,br)}}\in \overline{B}(z_0,br)^c\big)\\
&\le & c_1 \Phi(2r) \frac{c_2}{\Phi(br)} \E_x \tau_{D\cap B(z_0,br)} \le c_3 \E_x \tau_{D\cap B(z_0,br)}\,  .
\end{eqnarray*}
This finishes the proof. \qed

\noindent
{\bf Proof of Theorem \ref{t:main-finite}} Let $a\in (1/2,1)$ and choose $b:=\frac{2a}{1+2a}$ so that $b\in (1/2,a)$.
Let $D\subset B(z_0,r)$ and let $u$ be a non-negative function on $\X$ which is regular harmonic
with respect to $X$ in $D$ and vanishes on $B(z_0,r)\cap \left(\overline{D}^c\cup D^{\mathrm{reg}} \right)$.
Since $B(z_0, r/8)\subset B(z_0, r/4)\cap B(z_0, abr/2)$, it follows from Proposition \ref{p:main-fin-aux} and Lemma \ref{l:exit-comparison} that
$$
u(x)\asymp \E_x\tau_D \int_{\overline{B}(z_0,ar/2)^c} j(z_0,y)u(y)m(dy)\, ,\qquad x\in  D\cap B(z_0, r/8)\,  ,
$$
with a constant depending on $a$.
\qed

\bigskip
\noindent
{\bf Acknowledgements:} The main part of this paper was done during the visit of Panki Kim
to the University of Zagreb from June 28 to July 5, 2015. He thanks the Department of Mathematics
of the University of Zagreb for the hospitality.

\end{doublespace}

\bigskip
\noindent
\vspace{.1in}
\begin{singlespace}


\end{singlespace}

\vskip 0.1truein

\parindent=0em

{\bf Panki Kim}

Department of Mathematical Sciences and Research Institute of Mathematics,

Seoul National University, Building 27, 1 Gwanak-ro, Gwanak-gu Seoul 08826, Republic of Korea

E-mail: \texttt{pkim@snu.ac.kr}

\bigskip

{\bf Renming Song}

Department of Mathematics, University of Illinois, Urbana, IL 61801,
USA

E-mail: \texttt{rsong@math.uiuc.edu}

\bigskip

{\bf Zoran Vondra\v{c}ek}

Department of Mathematics, University of Zagreb, Zagreb, Croatia, and \\
Department of Mathematics, University of Illinois, Urbana, IL 61801,
USA

Email: \texttt{vondra@math.hr}
\end{document}